\documentclass[hidelinks,onefignum,onetabnum]{siamart251216}

\usepackage{lipsum}
\usepackage{amsfonts}
\usepackage{graphicx}
\usepackage{epstopdf}
\usepackage{algorithmic}
\usepackage{enumitem}
\usepackage{bm}

\ifpdf
  \DeclareGraphicsExtensions{.eps,.pdf,.png,.jpg}
\else
  \DeclareGraphicsExtensions{.eps}
\fi

\newsiamthm{assumption}{Assumption}
\newsiamremark{example}{Example}
\newsiamremark{remark}{Remark}
\newsiamremark{hypothesis}{Hypothesis}
\crefname{hypothesis}{Hypothesis}{Hypotheses}
\newsiamthm{claim}{Claim}
\newsiamremark{fact}{Fact}
\crefname{fact}{Fact}{Facts}

\headers{Convergence of primal-dual algorithm for GNE seeking}{W. Li, L. Pavel, P. J. Antsaklis\and H. Lin}

\title{On the convergence of discrete-time primal-dual algorithm for distributed GNE seeking\thanks{Submitted to the editors DATE.
}}

\author{Weijian Li\thanks{Department of Electrical Engineering, College of Engineering, University of Notre Dame, Notre Dame, IN, 46556, USA (\email{wli26@nd.edu}, \email{pantsakl@nd.edu}, \email{hlin1@nd.edu}).}
\and Lacra Pavel\thanks{ Department of Electrical and Computer Engineering, University of Toronto, Toronto, ON M5S 3G4, Canada (\email{pavel@control.toronto.edu}).}
\and Panos J. Antsaklis\footnotemark[2]
\and Hai Lin\footnotemark[2]
}

\usepackage{amsopn}

\ifpdf
\hypersetup{
pdftitle={On the convergence of discrete-time primal-dual algorithm for distributed GNE seeking},
pdfauthor={W. Li, L. Pavel, P. J. Antsaklis, H. Lin}
}
\fi

\begin{document}

\maketitle

% REQUIRED
\begin{abstract}
We investigate the convergence of a discrete-time distributed primal-dual algorithm for generalized Nash equilibrium (GNE) seeking in a class of noncooperative games with nonlinear coupling constraints. We establish semi-global convergence of the algorithm under a fixed stepsize, in the sense that convergence can be guaranteed from any prescribed compact set of initial conditions by choosing a sufficiently small stepsize.
Under a diminishing stepsize sequence, we further establish global convergence when the constraints are globally Lipschitz continuous. For the special case of affine constraints, global convergence is also guaranteed under a sufficiently small fixed stepsize. Finally, we extend these convergence results to a fully distributed GNE seeking algorithm in a partial-decision information setting, where each agent only has access to the actions of its neighbors. These results complement existing convergence results for continuous-time primal-dual dynamics, and provide new convergence guarantees for the discrete-time counterpart in distributed GNE seeking.
\end{abstract}

% REQUIRED

% REQUIRED

%%%%%%
%%%%%%--------------------
\section{Introduction}

Noncooperative games have been employed to model distributed decision-making in networked systems, with applications including demand-side management in smart grids \cite{saad2012game},  intelligent traffic control \cite{johansson2021strategic}, wireless communication networks \cite{yin2011nash}, economic systems \cite{yan2024incorporation} and social networks \cite{ghaderi2014opinion}.
In such games, each player is endowed with a local cost that depends on its own action as well as on the strategies of its opponents, while its feasible decision set may also be affected by the decisions of other agents when they compete for limited network resources.
A GNE is a reasonable solution since no player can decrease its local cost by unilaterally
changing its own decision \cite{facchinei2003finite, facchinei2010generalized}.
Thus, a variety of distributed algorithms have been developed for GNE seeking, including  primal-dual dynamics \cite{bianchi2021continuous, hu2022distributed}, operator splitting approaches \cite{yi2019operator, pavel2019distributed, pavel2026operator}, and pay-off based learning \cite{gao2022continuous, tatarenko2018learning}.

One of the most popular approaches for constrained optimization is the primal-dual algorithm since it is powerful to handle constraints, and can be easily implemented in a distributed manner. The algorithm was first introduced in \cite{arrow1958studies}, and its local convergence to a saddle point of the Lagrangian was subsequently established in \cite{polyak1970iterative}.
Then, global convergence was obtained for a modified version with perturbations and adaptive stepsizes \cite{kallio1999large}.
In \cite{nedic2009subgradient}, the primal-dual algorithm was applied to saddle-point problems, with global convergence guarantees under a diminishing stepsize sequence.
Following that, continuous-time primal-dual dynamics have also been developed for distributed optimization, and their asymptotic convergence has been extensively studied  \cite{wang2010control, cherukuri2017saddle, zeng2018distributed, kia2015distributed}.
The discrete-time counterparts were investigated in \cite{zhu2011distributed, lei2016primal}, where global convergence was established either under a diminishing stepsize sequence or for problems with local affine constraints.
It was further shown in \cite{bin2024semiglobal} that the discrete-time primal-dual algorithm for convex optimization with nonlinear constraints exhibits semi-global exponential convergence under a fixed stepsize, and the analysis was later extended to distributed optimization in \cite{ren2026exponentially}.

More recently, primal-dual algorithms have been developed for distributed Nash equilibria (NE) or GNE seeking. For instance, a continuous algorithm was proposed in \cite{gharesifard2013distributed} for subnetwork zero-sum games, while a discrete-time version was explored with switching communications in \cite{lou2015nash}.
In \cite{zeng2019generalized},  a distributed primal-dual algorithm was proposed for GNE seeking in multi-cluster games with nonsmooth payoff functions, nonlinear coupling constraints, and local set constraints.
In \cite{bianchi2021continuous}, a fully distributed continuous-time algorithm was proposed in a partial-decision information setting, and was further generalized to multi-integrator dynamic agents. 
Following that, primal-dual dynamics were characterized from a passivity-based perspective, and a class of generalized high-order dynamics was developed under relaxed monotonicity conditions in \cite{li2024passivity, li2025passivity}.
In practice, however, discretization is essential for the implementation of continuous-time dynamics. While the continuous-time primal-dual dynamics have been extensively investigated for distributed GNE seeking, the convergence of their discrete-time counterparts has not been well understood.

Motivated by the above observations, we investigate the convergence of a well-known discrete-time primal-dual algorithm for distributed GNE seeking. Our main contributions are summarized as follows.
\begin{enumerate}[label=\arabic*),  leftmargin=1.4em, itemsep=0.3em]
\item We establish semi-global convergence for this discrete-time distributed algorithm under a fixed stepsize. Specifically, for any prescribed compact set of initial conditions, there always exists a sufficiently small fixed stepsize such that the convergence is guaranteed. This result extends the existing analysis in \cite{bin2024semiglobal,ren2026exponentially} from constrained optimization to distributed GNE seeking with coupling constraints. Moreover, it complements the continuous-time convergence results in \cite{bianchi2021continuous,li2025passivity}.

\item We show global convergence of this algorithm under a diminishing stepsize sequence when the constraints are globally Lipschitz continuous. This result extends existing convergence guarantees for primal-dual algorithms in saddle-point problems and two-player zero-sum games to distributed GNE seeking in multi-player noncooperative games with nonlinear coupling constraints \cite{nedic2009subgradient,zhu2011distributed,lou2015nash}.

\item For games with affine constraints, we further establish global convergence of the discrete-time algorithm under a fixed stepsize. This strengthens the semi-global convergence result since the stepsize is independent of the initial conditions.
The algorithm differs from the operator-splitting approaches in \cite{yi2019operator,pavel2019distributed,belgioioso2022distributed,pavel2026operator} in that it does not require correction steps, thereby facilitating distributed implementation.
To our knowledge, such a convergence result is obtained for the first time.

\item By establishing the restricted strong monotonicity of an augmented pseudogradient mapping, we generalize the above results to a partial-decision information setting, where each agent only has access to the actions of its neighbors.
\end{enumerate}

The rest of this paper is organized as follows. In Section \ref{sec:preliminary}, we review the mathematical preliminaries. In Section \ref{sec:setup}, we introduce a discrete-time primal-dual algorithm for distributed GNE seeking and present the problem statements.
Following that, in Sections \ref{sec:semi} and \ref{sec:global}, we analyze the semi-global convergence under a fixed stepsize and the global convergence under a diminishing stepsize. In Section \ref{sec:affine}, we further show that, in the case of affine constraints, the algorithm is globally convergent under a fixed stepsize. In Section \ref{sec:partial}, we extend the results to partial-decision information settings. We provide simulations for illustration in Section \ref{sec:numerical}. Finally, we close this paper with concluding remarks in Section \ref{sec:conclusion}.

%%%%%%
%%%%%%--------------------
\section{Terminologies and basics}
\label{sec:preliminary}

In this section, we introduce relevant notation and review necessary concepts on projections, monotone operators, and graph theory.

\emph{Notation:}
Let $\mathbb R_+^m$ be the set of $m$-dimensional nonnegative real column vectors, 
and $\mathbb N$ be the set of nonnegative integers.
Denote by $1_m (0_m)$ the $m$-dimensional vector with all entries of $1 (0)$, and by $I_n$ the $n$-by-$n$ identity matrix. 
We simply write $\bm 0$ for zero vectors/matrices of appropriate dimensions when there is no confusion.
Let $(\cdot)^\top$, $\otimes$, and $\|\cdot\|$ denote the transpose, the Kronecker product, and the Euclidean norm of a vector or the induced $\ell_2$ norm of a matrix.
Denote by ${\rm Ker}(\cdot)$ and ${\rm rank}(\cdot)$ the kernel space and the rank of a matrix.
Let $\mathcal I = \{1, \dots, N\}$,  ${\rm col}\{x_i\}_{i \in \mathcal I} := [x_1^\top, \dots, x_N^\top]^\top \in \mathbb R^n$ for $x_i \in \mathbb R^{n_i}$ with $n = \sum_{i \in \mathcal I} n_i$, and  $X := {\rm blkdiag}\{X_i\}_{i \in \mathcal I}$ be a block diagonal matrix with blocks $X_i \in \mathbb R^{p_i \times q_i}$.
For a set $\mathcal Q \subset \mathbb R^n$, $\bar{\mathcal Q}$ denotes its closure.
For $x_0 \in \mathbb R^n$ and $K > 0$, we define $\bar{\mathbb B}_K(x_0) := \{x \in \mathbb R^n \mid \|x - x_0\| \le K\}$. For a differentiable function $J(x, y)$, $\nabla_x J(x, y)$ is the partial gradient of $J$ with respect to $x$.

\emph{Projection:} 
Let $\Omega \subset \mathbb R^n$ be a nonempty closed convex set such that $\lambda x + (1 - \lambda)y \in \Omega$ for all $x, y \in \Omega$ and $\lambda \in [0, 1]$.
The normal cone to $\Omega$ at $x \in \Omega$ is defined as $\mathcal N_{\Omega}(x) := \{v \in \mathbb R^n \mid v^\top(y - x) \le 0, \forall y \in \Omega\}$.
For $z \in \mathbb R^n$, the Euclidean projection of $z$ onto $\Omega$ is given by
$P_\Omega(z) := {\rm argmin}_{y \in \Omega} \|y - z\|$.
Referring to \cite{facchinei2003finite}, the projection operator $P_\Omega$ is firmly nonexpansive, i.e.,
\begin{subequations}
\begin{align}
&\|P_\Omega(z) - P_{\Omega}(x)\| \le \|z - x\|,~\forall z, x \in \mathbb R^n, \label{projection1}\\
&\|P_\Omega(z) - P_{\Omega}(x)\|^2 \le \big\langle P_\Omega(z) - P_{\Omega}(x),z - x\big\rangle,~\forall z, x \in \mathbb R^n. \label{projection2}
\end{align}
\end{subequations}

\emph{Monotone operator:}
An operator (or mapping) $F: \mathbb R^n \rightarrow \mathbb R^n$ is monotone  if $\langle x - y, F(x) - F(y)\rangle \ge 0$, and $\mu$-strongly monotone if $\langle x - y, F(x) - F(y)\rangle \ge \mu \Vert x - y \Vert^2$, for all $x, y \in \mathbb R^n$.
It is $\theta$-Lipschitz continuous if $\Vert F(x) - F(y)\Vert \le \theta \Vert x - y \Vert, \forall x, y \in \mathbb R^n$. Let $f: \mathbb R^n \rightarrow \mathbb R$ be a continuously differentiable convex function and $\nabla f(x)$ be the gradient of $f$ at $x$. Then
$f(y) \ge f(x) + \langle y - x, \nabla f(x)\rangle, \forall x, y \in \mathbb R^n$, and  $\nabla f$ is monotone. 

\emph{Graph theory:}
Consider a multi-agent network modeled by an undirected graph $\mathcal G(\mathcal I, \mathcal E, \mathcal A)$, where $\mathcal I$ is the node set, $\mathcal E \subset \mathcal I \times \mathcal I$ is the edge set, and $\mathcal A = [a_{ij}]\in \mathbb R^{N \times N}$ is the adjacency matrix such that $a_{ij} = a_{ji} > 0$ if nodes $i$ and $j$ are connected, and $a_{ij} = 0$ otherwise. Suppose that there are no self-loops in $\mathcal G$, i.e., $a_{ii} = 0$.
Node $j$ is a neighbor of node $i$ if and only if $a_{ij} > 0$.
Let $\mathcal N_i$ be the set of node $i$’s neighbors. 
The Laplacian matrix $L_G$ is $L_G = \mathcal D - \mathcal A$, where $\mathcal D = {\rm diag}\{d_i\}\in \mathbb R^{N \times N}$, and $d_i = \sum_{j \in \mathcal I} a_{ij}$.
Graph $\mathcal G$ is connected if there exists a path between every pair of distinct nodes. 
If $\mathcal G$ is connected, then $L_G = L_G^\top$, ${\rm rank}(L_G) = N - 1$, and ${\rm Ker}(L_G) = \{k 1_N: k \in \mathbb R\}$.

%%%%%%
%%%%%%--------------------
\section{Formulation and algorithm}
\label{sec:setup}

In this section, we formulate the distributed GNE seeking problem, introduce the discrete-time primal-dual algorithm under a full-decision information setting, and give the problem statements.

%%%%%%
%%%%%%--------------------
\subsection{Game setup}

Consider a set $\mathcal I = \{1, \dots, N\}$ of $N$ players (agents) participating in a noncooperative game. Player $i$ controls its action (strategy) $x^i \in \mathbb R^{n_i}$.
Let $x = (x^i, x^{-i}) \in \mathbb R^n$ denote the joint action profile of all players, where 
$\sum_{i = 1}^N n_i = n$, and $x^{-i}$ is the $(N-1)$-tuple of all agents' actions except agent $i$'s decision. Alternatively, $x = {\rm col}\{x^i\}_{i \in \mathcal I} \in \mathbb R^n$.
Player $i$ aims to minimize its local cost $J_i(x^i, x^{-i}): \mathbb R^n \rightarrow \mathbb R$, which depends on its own strategy $x^i$ and the actions of its opponents $x^{-i}$.
In addition, actions of the players are coupled through separable constraints.
Specifically, the feasible set is given by
$X := \{x \in \mathbb R^n \mid g(x) \le 0_m\}$, where $g(x) := \sum_{i\in \mathcal I} g_i(x^i)$, and $g_i: \mathbb R^{n_i} \rightarrow \mathbb R^m$ is a private function known only to agent $i$.
Given $x^{-i}$, the feasible decision set of player $i$ is 
$X_i(x^{-i}) := \{x^i \in \mathbb R^{n_i}\mid (x^i, x^{-i}) \in X\}$.
Accordingly, player $i$ seeks a best response strategy by solving
\begin{equation}
\label{formulation}
\min_{x^i\in \mathbb R^{n_i}} J_i(x^i, x^{-i}),~
{\rm s.t.}~x^i \in X_i(x^{-i}).
\end{equation}

A collective strategy $x^* = (x^{i, *}, x^{-i,*})$ is called a GNE if 
\begin{equation*}
\label{def:GNE}
x^{i, *} \in {\rm argmin}_{x^i} J_i(x^i, x^{-i,*}),
~{\rm s.t.}~(x^i, x^{-i,*}) \in X,~\forall i \in \mathcal I.
\end{equation*}
In the absence of coupled constraints, $x^*$ is an NE.

To ensure the well-posedness of this game, we make the following standing assumptions \cite{pavel2019distributed, yi2019operator, bianchi2021continuous}.

\begin{assumption}
\label{ass:convex}
\leavevmode
\begin{enumerate}[
label=\textnormal{(A\arabic*)},
leftmargin=*,
labelindent=\parindent,
align=left
]
\item For $i \in \mathcal I$, $J_i$ is continuously differentiable and convex in $x^i$, given $x^{-i}$.

\item For $i \in \mathcal I$, $g_i$ is twice continuously differentiable and convex.

\item The set $X$ is nonempty and satisfies Slater’s constraint qualification.
\end{enumerate}
\end{assumption}

The Lagrangian of agent $i$ is defined as
\begin{equation}
\label{largrangian}
\mathcal L_i(x^i, \lambda^i; x^{-i}) := J_i(x^i, x^{-i}) + \lambda^{i, \top} g(x),
\end{equation}
with dual variable (multiplier) $\lambda^i \in \mathbb R^m$.
When $x^{i,*}$ is an optimal solution to \eqref{formulation}, there is $\lambda^{i,*}$ such that the following Karush-Kuhn-Tucker (KKT) conditions hold \cite{facchinei2010generalized}:
\begin{equation}
\label{GNE}
\begin{aligned}
&0_{n_i} = \nabla_{x^i} \mathcal L_i(x^{i,*}, \lambda^{i,*}; x^{-i,*}), ~\forall i \in \mathcal I, \\
&g(x^*) \le 0, ~\lambda^{i, *} \ge 0,~{\rm and}~\big\langle\lambda^{i, *}, g(x^*) \big\rangle = 0.
\end{aligned}
\end{equation}

In general, the players may have different multipliers at a GNE, i.e., $\lambda^{i,*} \neq\lambda^{j,*}$ for some $i\neq j$. In this paper, we focus on seeking a GNE characterized by the same
Lagrangian multiplier, named variational GNE (v-GNE), which has the economic interpretation of no price discrimination \cite{kulkarni2012variational}.
A v-GNE $x^* \in X$ of the game \eqref{formulation} is a solution to the following variational inequality:
\begin{equation*}
{\rm VI}(F, X): \big\langle x - x^*, F(x^*)\big\rangle \ge 0, ~\forall x \in X,
\end{equation*}
where $F$ denotes the pseudo-gradient mapping defined by
\begin{equation}
\label{pseudo-gradient}
F(x) := {\rm col}\{\nabla_{x^i} J_i(x^i, x^{-i})\}_{i \in \mathcal I} \in \mathbb R^n.
\end{equation}

Note that $x^*$ solves ${\rm VI}(F, X)$ if and only if there exists a multiplier $\lambda_c^*$ such that the following KKT conditions are satisfied \cite{facchinei2010generalized}:
\begin{equation}
\label{vGNE}
\begin{aligned}
&0_{n_i} = \nabla_{x^i} J_i(x^{i,*}, x^{-i,*}) + \nabla g_i(x^{i, *})^\top \lambda_c^*, ~\forall i \in \mathcal I, \\
&0_m \in - g(x^*) + \mathcal N_{\mathbb R^m_+}(\lambda_c^*),
\end{aligned}
\end{equation}
where $\mathcal N_{\mathbb R^m_+}(\lambda_c^*)$ is the normal cone of $\mathbb R^m_+$ at $\lambda_c^*$.
We impose the following assumption on the pseudogradient $F$ \cite{pavel2019distributed, yi2019operator}.
\begin{assumption}
\label{ass:monotone}
The mapping $F$ is $\mu$-strongly monotone and $\theta$-Lipschitz continuous.
\end{assumption}

%%%%%%
%%%%%%--------------------
\subsection{Distributed primal-dual algorithm}

We suppose that each agent $i$ has access only to its local private information, i.e., $J_i$ and $g_i$, but knows the actions $x^{-i}$ of all its opponents, which corresponds to a full-decision information setting \cite{hu2022distributed, yi2019operator, zeng2019generalized}.
Then, agent $i$ can compute the exact partial-gradient $\nabla_{x^i} J_i(x^i, x^{-i})$.

We now introduce the necessary notation for the algorithm design.
Agent $i$ controls its local strategy $x^i \in \mathbb R^{n_i}$, and keeps a local copy of the multiplier $\lambda^i \in \mathbb R^m$ as an estimate of $\lambda_c^*$. Besides, it maintains a local auxiliary variable $z^i \in \mathbb R^m$ to facilitate consensus among the multipliers $\lambda^i$.
The communication topology among the agents is modeled by an undirected graph $\mathcal G_c$, where agent $i$ exchanges its local variables $\{\lambda^i, z^i\}$ with its neighbors $j \in \mathcal N_i$.
We make the following assumption on $\mathcal G_c$.

\begin{assumption}
\label{ass:graph}
The undirected graph $\mathcal G_c$ is connected.
\end{assumption}

With the above preparations, we are ready to present the distributed discrete-time primal-dual algorithm for v-GNE seeking as
\begin{equation}
\label{alg:primal-dual:agent}
\left\{
\begin{aligned}
x^i_{t+1} &=  x^i_t - \gamma \nabla_{x^i} J_i(x_t^i, x_t^{-i}) - \gamma\nabla g_i(x_t^i)^\top \lambda_t^i,\\ 
z^i_{t+1} &= z^i_t + \gamma \sum_{j \in \mathcal N_i} a_{ij}(\lambda_t^i - \lambda_t^j),  \\ 
\lambda^i_{t+1} &= P_{\mathbb R^m_+}\big[\lambda^i_t + \gamma g_i(x_t^i) -  \gamma \sum_{j \in \mathcal N_i} a_{ij}(z_t^i - z_t^j) - \gamma \sum_{j \in \mathcal N_i} a_{ij}(\lambda_t^i - \lambda_t^j)\big],
\end{aligned}
\right.
\end{equation}
where  $t \in \mathbb N$ denotes the iteration index, $\gamma > 0$ is the stepsize, $P_{\mathbb R^m_+}(\cdot)$ denotes the Euclidean projection onto $\mathbb R^m_+$, and $a_{ij}$ is the $(i, j)$-th entry of the adjacency matrix of $\mathcal G_c$. The initial conditions satisfy $x^i_0 \in \mathbb R^{n_i}$, $z^i_0 \in \mathbb R^{m}$, and $\lambda^i_0 \in \mathbb R^m_+$.

Let $\tilde m = Nm$ and $L_G$ be the Laplacian matrix of $\mathcal G_c$. We further define
\begin{equation*}
\begin{aligned}
G(x) &:= {\rm col} \{g_i(x^i)\}_{i \in \mathcal I} \in \mathbb R^{\tilde m},
~\nabla G(x) := {\rm blkdiag}\{\nabla g_i(x^i)\}_{i \in \mathcal I}  \in \mathbb R^{\tilde m \times n}, \\
\lambda &:= {\rm col} \{\lambda^i\}_{i \in \mathcal I} \in \mathbb R^{\tilde m},
~z := {\rm col} \{z^i\}_{i \in \mathcal I} \in \mathbb R^{\tilde m},
~L := L_G \otimes I_m \in \mathbb R^{\tilde m \times \tilde m}.
\end{aligned}
\end{equation*}
Then, in a compact form, algorithm \eqref{alg:primal-dual:agent} reads as
\begin{equation}
\label{alg:primal-dual}
\left\{
\begin{aligned}
x_{t+1} &= x_t - \gamma F(x_t) - \gamma \nabla G(x_t)^\top \lambda_t, \\
z_{t+1} &= z_t + \gamma L \lambda_t, \\
\lambda_{t+1} &= P_{\mathbb R^{\tilde m}_+} \big[\lambda_t + \gamma G(x_t) - \gamma L z_t - \gamma L \lambda_t \big],
\end{aligned}
\right.
\end{equation}
where $x_0 \in \mathbb R^n$, $z_0 \in \mathbb R^{\tilde m}$, and $\lambda_0 \in \mathbb R^{\tilde m}_+$.

Intuitively, agent $i$ updates $x^i$ by performing gradient descent on its local Lagrangian $\mathcal L_i$ in \eqref{largrangian}, and updates $\lambda^i$ via a gradient ascent step combined with a proportional-integral synchronization mechanism involving $z^i$.
The following lemma clarifies the relationship between the v-GNE of \eqref{formulation} and the equilibria of \eqref{alg:primal-dual}. Its proof follows similar arguments to those in \cite[Lemma 2]{bianchi2021continuous} and \cite[Theorem 4.1]{zeng2019generalized}, and is therefore omitted.

\begin{lemma}
\label{lem:KKT-EQ}
Under Assumptions \ref{ass:convex}, \ref{ass:monotone} and \ref{ass:graph}, a profile $x^*$ is a v-GNE of (\ref{formulation}) if and only if there exists $(\lambda^*, z^*)$ such that $(x^*, z^*, \lambda^*)$ is an equilibrium of \eqref{alg:primal-dual}.
\end{lemma}

\begin{remark}
The continuous-time counterpart of \eqref{alg:primal-dual} was proposed in \cite{bianchi2021continuous}, where its asymptotic convergence was established. Subsequently, this dynamics was generalized within a passivity-based framework in \cite{li2025passivity}. While discrete-time algorithms are essential for numerical implementation, the convergence of such a discrete-time primal-dual algorithm for the game \eqref{formulation} has not been investigated, to the best of our knowledge.
In general, the global convergence of \eqref{alg:primal-dual} with a fixed stepsize is not guaranteed, as illustrated by the following example.
\end{remark}

\begin{example}
\label{Ex1}
Consider a network with two agents, where
$J_i(x^i, x^{-i}) = (x^i)^2 + + \rho x^i x^j$, and $g_i = (x^i)^2 - 1$, $x^i \in \mathbb R$, and $0 < \rho < 2$. The v-GNE is $x^{1, *} = x^{2,*} = 0$.
Take the initial point $x^1_0 = x^2_0 = q_0$, $z_0^1 = z_0^2 = 0$, and $\lambda_0^1 = \lambda_0^2 = \eta_0$. It follows from \eqref{alg:primal-dual:agent} that $x^1_t = x^2_t =: q_t$,
$\lambda^1_t = \lambda^2_t =: \eta_t$, and $z_t^1 = z_t^2 = 0$.
Moreover, we have
\begin{equation}
\label{alg:Ex1}
\left\{
\begin{aligned}
q_0 & \ge \sqrt{2}, \\
\eta_0 &\ge  \frac{1 + \sqrt{2}}{2 \gamma}
\end{aligned}
\right.
\implies
\left\{
\begin{aligned}
q^2_{t+1} &= (1 - 2 \gamma - \rho \gamma - 2 \gamma \eta_t)^2 q_t^2 \ge 2 q_t^2 \ge 2, \\
\eta_{t+1} &= \max\{0, \eta_t + \gamma(q_t^2 - 1)\} \ge \eta_t \ge \frac{1 + \sqrt{2}}{2 \gamma}. 
\end{aligned}
\right.
\end{equation}
Thus, an arbitrarily small fixed stepsize does not, in general, guarantee the convergence of \eqref{alg:primal-dual} to the v-GNE.
\end{example}

The above observation motivates a further investigation into the convergence of \eqref{alg:primal-dual}. In particular, we are interested in the following questions:

\begin{enumerate}[label=\roman*)]
\item What convergence guarantees can be established for \eqref{alg:primal-dual} with a fixed stepsize?

\item Under what conditions can global convergence be achieved?

\item Are there special cases for which global convergence can be guaranteed?

\item Is it possible to extend the convergence results to a partial-decision information setting?
\end{enumerate}

%%%%%%
%%%%%%--------------------
\section{Semi-global convergence}
\label{sec:semi}

In this section, we establish the semi-global convergence of algorithm \eqref{alg:primal-dual} under a fixed stepsize. We first present the main result, and then provide its proof.

\begin{theorem}
\label{thm:semi_global} 
Suppose that Assumptions \ref{ass:convex}, \ref{ass:monotone} and \ref{ass:graph} hold. Then, for every compact set $\Xi_0 \subset \mathbb R^n \times \mathbb R^{\tilde m} \times \mathbb R_+^{\tilde m}$ of initial conditions, there exists $\bar \gamma^s > 0$ such that, for any fixed stepsize $\gamma \in (0, \bar \gamma^s)$ and any initial condition $(x_0, z_0, \lambda_0) \in \Xi_0$, the sequence $\{x_t\}$ generated by \eqref{alg:primal-dual} converges to the v-GNE $x^*$ of the game \eqref{formulation}.
\end{theorem}

\begin{remark}
Theorem \ref{thm:semi_global} establishes a semi-global convergence result for \eqref{alg:primal-dual}.
In particular, for an arbitrarily large compact set $\Xi_0$, there always exists a constant ${\bar \gamma}^s$ such that the convergence is guaranteed for all initial conditions in $\Xi_0$ whenever $\gamma \in (0, {\bar \gamma}^s)$. In general, ${\bar \gamma}^s$ depends on $\Xi_0$, and a smaller stepsize may be required as $\Xi_0$ becomes larger. This differs from global convergence, for which the stepsize can be chosen independent of the initial conditions \cite{grne2013nonlinear}. 
\end{remark}

\begin{remark}
The semi-global stability of a discrete-time primal-dual algorithm for constrained optimization was established in \cite{bin2024semiglobal} and subsequently extended to distributed optimization in \cite{ren2026exponentially}. 
Our work differs from these studies in two aspects. On the one hand, we consider a noncooperative game setting, where the cost of each player depends not only on its own action but also on the actions of others. On the other hand, we consider coupled constraints in \eqref{formulation}, which require the distributed design to enforce the shared constraints and reach consensus on the dual variables, and lead to a more involved convergence analysis.
To the best of our knowledge,  the semi-global convergence of \eqref{alg:primal-dual} has not appeared in any existing literature.
\end{remark}

%%%%%%
%%%%%%--------------------
\subsection{Technical preliminaries}

To prove Theorem \ref{thm:semi_global}, we introduce some preliminary definitions and construct a Lyapunov function candidate.

Let $(x^*, z^*, \lambda^*)$ be an equilibrium of \eqref{alg:primal-dual}, and define the error variables as
\begin{equation}
\label{error:var}
\tilde x: = x - x^*, ~\tilde z := z - z^*,~{\rm and}~\tilde \lambda := \lambda - \lambda^*.
\end{equation}

We fix an arbitrary compact set $\Xi_0 \subset \mathbb R^n \times \mathbb R^{\tilde m} \times \mathbb R_+^{\tilde m}$. 
Then there exists $K_0 > 0$ such that $\Xi_0 \subset \bar{\mathbb B}_{K_0}(x^*, z^*,  \lambda^*)$.
We further choose $K \ge 2K_0 + 1$, and define
$$
\mathcal K := \bar{\mathbb B}_K(x^*, z^*, \lambda^*),~{\rm and}~
\mathcal K_x := \{x \in \mathbb R^n: \exists (z, \lambda) \in \mathbb R^{\tilde m} \times \mathbb R_+^{\tilde m}, (x, z, \lambda) \in \mathcal K\}.
$$

For the ease of notation, we denote by
\begin{equation}
\begin{aligned}
\label{def:Delta}
\Delta_1(x, \lambda) := F(x) + \nabla G(x)^\top \lambda, ~{\rm and}~
\Delta_2(x, z, \lambda) := G(x) - L z - L\lambda.
\end{aligned}
\end{equation}
Note that $\Delta_1(x^*, \lambda^*) = F(x^*) + \nabla G(x^*)^\top \lambda^* = \mathbf{0}$. Since $\mathcal K$ is compact, under Assumptions \ref{ass:convex} and \ref{ass:monotone},  there exist scalars $k_1, \dots k_6 > 0$ such that
\begin{equation}
\begin{aligned}
\label{def:k:1:6}
&k_1 > \sup_{(x, z, \lambda) \in \mathcal K} \| \Delta_1(x, \lambda) \|,~
k_2 > \sup_{x \in \mathcal K_x} \|\nabla G(x) \|, ~
k_3 > \sup_{x \in \mathcal K_x} \|G(x)\|,\\
&k_4 > \sup_{(x,z,\lambda) \in \mathcal K} \|\Delta_2(x,z,\lambda)-\Delta_2(x^*,z^*,\lambda^*)\|,\\
&\forall x \in \mathcal K_x, ~\|\Delta_1(x, \lambda^*)\| \le k_5\|\tilde x\|,\\
&\forall x, y \in \mathcal K_x, ~\|\nabla G(x) - \nabla G(y)\| \le k_6\|x - y\|.
\end{aligned}
\end{equation}
Let $k_7 = \max\{k_5, 1\}$ and $k_L := \|L\|$.

Finally, we fix a constant $\kappa_s > 1$, and choose the stepsize $\gamma$ such that
\begin{equation}
\label{gamma:s}
\gamma \in (0, \bar \gamma^s), ~\bar \gamma^s:= \min_{j \in \{1, \dots, 5\}}\{\bar \gamma_j^s\},
\end{equation}
where
\begin{subequations}
\begin{align}
\bar \gamma_1^s &:= \min\Bigg\{\frac{1}{24K_0k_1}, \frac{1}{\sqrt{6}k_1}, 
\frac{1}{48K_0^2 k_L}, \frac{1}{2\sqrt{6}K_0k_L}, \frac{1}{24K_0k_4}, \frac{1}{\sqrt{6}k_4}\Bigg\}, \label{gamma:s:1}\\
\bar \gamma_2^s &:= \min\Bigg\{\frac{1}{4 \kappa_s k_2}, \frac{1}{4 k_L}\Bigg\},
\label{gamma:s:2}\\
\bar{\gamma}_3^s &:= \min\Bigg\{\frac{\mu}{3\big(8(1+\kappa_s)^2k_2^2 + 2\kappa_s k_5 k_6K +\kappa_s^2 k_6^2 K^2 +\kappa_s k_5^2\big)},
\label{gamma:s:3}\\
&\hspace{2.8cm}
\sqrt{\frac{\mu}{
9\kappa_s k_6 k_7 K (k_5 + 1)}},
\sqrt[3]{\frac{\mu}{
48\kappa_s^2k_2^2k_7^2}}\Bigg\}, \notag\\
\bar{\gamma}_4^s &:= \min\Bigg\{ \frac{\kappa_s -1} {12\kappa_s k_6 k_7(k_5 + 1)K},
\frac{\sqrt{\kappa_s -1}} {8\kappa_s k_2 k_7} \Bigg\}, \label{gamma:s:4}\\
\bar{\gamma}_5^s &:= \min\Bigg\{\frac{1}{22 k_L}, \frac{1}{\sqrt[3]{16}\, k_L}\Bigg\}. \label{gamma:s:5}
\end{align}
\end{subequations}

\begin{remark}
The admissible stepsize $\gamma$ in \eqref{gamma:s} depends on the initial conditions, properties of the pseudogradient mapping and the coupling constraints, and the communication topology. In particular, the constants involved in the stepsize bound are determined over a compact set, whose size depends on the prescribed set of initial conditions. 
This dependence explains the semi-global convergence guarantee in Theorem \ref{thm:semi_global}. We emphasize that the stepsize conditions in \eqref{gamma:s} are sufficient for convergence.
\end{remark}

The next lemma characterizes the behavior of \eqref{alg:primal-dual} after one-step update.
\begin{lemma}
\label{lem:Iter:Bound}
Let $\gamma$ satisfy \eqref{gamma:s}.
Under Assumptions \ref{ass:convex}, \ref{ass:monotone}, and \ref{ass:graph}, it holds that
\begin{equation}
\label{one-step}
(x_t, z_t, \lambda_t) \in \bar{\mathbb B}_{2K_0}(x^*, z^*, \lambda^*) \implies
(x_{t+1}, z_{t+1}, \lambda_{t+1}) \in \mathcal K.
\end{equation}
\end{lemma}

\begin{proof}
By \eqref{alg:primal-dual}, \eqref{def:k:1:6} and \eqref{gamma:s:1}, we obtain
\begin{equation*}
\begin{aligned}
\|\tilde x_{t+1}\|^2 &= \|\tilde x_t\|^2 - 2 \gamma \big\langle \tilde x_t, \Delta_1(x_t, \lambda_t) \big\rangle + \gamma^2 \|\Delta_1(x_t, \lambda_t)\|^2 \\
&\le \|\tilde x_t\|^2 + 4\gamma K_0 k_1 + \gamma^2 k_1^2 \le \|\tilde x_t\|^2 + 1/3.
\end{aligned}
\end{equation*}
Similarly, using $L \lambda^* = \mathbf{0}$, we have
\begin{equation*}
\begin{aligned}
\|\tilde z_{t+1}\|^2  
&\le \|\tilde z_t\|^2 + 8 \gamma K_0^2 k_L + 4 \gamma^2 K_0^2 k_L^2 \le \|\tilde z_t\|^2 + 1/3.
\end{aligned}
\end{equation*}
Since $\lambda^* = P_{\mathbb R^{\tilde m}_+} \big[\lambda^* + \gamma \Delta_2(x^*,z^*,\lambda^*) \big]$, recalling \eqref{projection1} gives
\begin{equation*}
\|\tilde \lambda_{t+1}\|^2 \le 
\|\tilde \lambda_t\|^2 + 4 \gamma K_0 k_4 + \gamma^2 k_4^2 \le \|\tilde \lambda_t\|^2 + 1/3.
\end{equation*}

Therefore, we derive
$$
\|\tilde x_{t+1}\|^2 + \|\tilde \lambda_{t+1}\|^2 + \|\tilde z_{t+1}\|^2 \le \|\tilde x_t\|^2 + \|\tilde \lambda_t\|^2 + \|\tilde z_t\|^2 + 1,
$$
which implies \eqref{one-step}. This completes the proof.
\end{proof}

Consider a Lyapunov candidate as
\begin{equation}
\label{lyap:V:s}
V_s(x, z, \lambda) = V_{s,1}(x)  + V_{s,2}(x, \lambda)  + V_{s,3}(z, \lambda),
\end{equation}
where 
\begin{equation*}
\begin{aligned}
V_{s,1}(x) &= \|\tilde x\|^2, \\
V_{s,2}(x, \lambda) &= 2 \gamma \kappa_s \big\langle \tilde x, \nabla G(x)^\top \tilde \lambda\big\rangle, \\
V_{s,3}(z, \lambda) &= \|\tilde z\|^2 + \|\tilde \lambda\|^2 + 2 \gamma \big\langle \tilde \lambda, L \tilde z \big\rangle.
\end{aligned}
\end{equation*}

\begin{remark}
The Lyapunov function in \eqref{lyap:V:s} is constructed by incorporating two cross terms into a standard quadratic Lyapunov candidate to account for the coupling among the primal, dual and auxiliary variables. The first cross term $V_{s,2}$ is introduced to provide additional negative terms in the Lyapunov difference \cite{bin2024semiglobal}, which are essential for establishing the convergence of \eqref{alg:primal-dual}.
The second cross term $2\gamma \langle \tilde \lambda, L \tilde z \rangle$ captures the interaction between the dual and auxiliary variables.
\end{remark}

The following lemma shows that $V_s$ is positive definite over the compact set $\mathcal K$. The proof is straightforward in view of \eqref{def:k:1:6} and \eqref{gamma:s:2}, together with the Cauchy-Schwarz and Young inequalities, and is therefore omitted.
\begin{lemma}
\label{lem:lyap:bound}
Let $\gamma$ satisfy \eqref{gamma:s}.
Under Assumptions \ref{ass:convex},  \ref{ass:monotone} and \ref{ass:graph}, it holds that
\begin{equation*}
\label{Lyap:s:PD}
\forall (x, z, \lambda) \in \mathcal K, ~\frac 12 \|(\tilde x, \tilde z, \tilde \lambda)\|^2 \le V_s(x, z, \lambda) 
\le \frac 32  \|(\tilde x, \tilde z, \tilde \lambda)\|^2.
\end{equation*}
\end{lemma}

%%%%%%
%%%%%%--------------------
\subsection{Convergence analysis}
\label{sec:semi:pf}
In this subsection, we complete the proof of Theorem \ref{thm:semi_global}. We first establish the following one-step descent property of $V_s$ in \eqref{lyap:V:s}.

\begin{lemma}
\label{lem:lyap:V:onestep}
Consider algorithm \eqref{alg:primal-dual}, and let the Lyapunov function $V_s$ be defined in \eqref{lyap:V:s}. Suppose that Assumptions \ref{ass:convex},  \ref{ass:monotone} and \ref{ass:graph} hold, and $\gamma$ satisfies \eqref{gamma:s}. If $(x_t, z_t, \lambda_t) \in \bar{\mathbb B}_{2K_0}$, then
\begin{equation}
\begin{aligned}
\label{lyap:Vs:onestep}
&V_s(x_{t+1}, z_{t+1}, \lambda_{t+1}) - V_s(x_t, z_t, \lambda_t) \\
&\le -\gamma \mu\|\tilde x_t\|^2  
- \frac 12 \|\tilde \lambda_{t+1}- \tilde \lambda_t\|^2 
- \frac{\gamma^2(\kappa_s - 1)}{2} \|\nabla G(x_t)^\top \tilde \lambda_t\|^2 - \frac{\gamma}{k_L}\|L \tilde \lambda_t\|^2.
\end{aligned}
\end{equation}
\end{lemma}

\begin{proof}
By Lemma \ref{lem:Iter:Bound}, $(x_{t+1}, z_{t+1}, \lambda_{t+1}) \in \mathcal K$.
Since $\Delta_1(x^*, \lambda^*) = \mathbf{0}$, we have
\begin{equation}
\begin{aligned}
\label{property:Vs1}
&V_{s,1}(x_{t+1}) - V_{s,1}(x_t) = -2 \gamma \big\langle \tilde x_t, \Delta_1(x_t, \lambda_t) \big\rangle 
+ \gamma^2 \Vert \Delta_1(x_t, \lambda_t)\Vert^2 \\
&\le -2 \gamma \mu \Vert \tilde x_t \Vert^2 
- 2\gamma \big\langle \tilde x_t, \nabla G(x_t)^\top \tilde \lambda_t \big\rangle - 2\gamma \big\langle \tilde x_t, \big(\nabla G(x_t) - \nabla G(x^*) \big)^\top \lambda^* \big\rangle \\
&~~~ + \gamma^2 \Vert \Delta_1(x_t, \lambda_t)\Vert^2,
\end{aligned}
\end{equation}
where the inequality follows from Assumption \ref{ass:monotone}.

Define $R_t := \big(\nabla G(x_{t+1}) - \nabla G(x_t)\big)^\top \tilde \lambda_{t+1}$.
It is then straightforward that
\begin{equation}
\begin{aligned}
\label{property:Vs2}
V_{s,2}&(x_{t+1}, \lambda_{t+1}) - V_{s,2}(x_t, \lambda_t) \\
=&~ 2\gamma \kappa_s \big\langle \tilde x_t, \nabla G(x_t)^\top (\tilde \lambda_{t+1} - \tilde \lambda_t) \big\rangle + 2 \gamma \kappa_s \big\langle \tilde x_t, R_t\big\rangle \\
&- 2 \gamma^2 \kappa_s \|\nabla G(x_t)^\top \tilde \lambda_t\|^2 
- 2 \gamma^2 \kappa_s \big\langle \Delta_1(x_t, \lambda^*), \nabla G(x_t)^\top\tilde{\lambda}_t  \big\rangle \\
&- 2 \gamma^2 \kappa_s \big\langle \Delta_1(x_t, \lambda_t), \nabla G(x_t)^\top(\tilde \lambda_{t+1} - \tilde \lambda_t) + R_t\big\rangle.
\end{aligned}
\end{equation}

Notice that $\lambda^* = P_{\mathbb R_+^{\tilde m}}[\lambda^* + \gamma\Delta_2(x^*, z^*, \lambda^*)]$ and 
$2 \big\langle \tilde \lambda_{t+1}, \tilde \lambda_t \big\rangle = \Vert \tilde \lambda_{t+1} \Vert^2 + \Vert \tilde \lambda_t \Vert^2 - \Vert \tilde \lambda_{t+1} - \tilde \lambda_t\Vert^2$.
Recalling \eqref{projection2} yields
\begin{equation}
\begin{aligned}
\label{lambda:s:descent}
\Vert \tilde \lambda_{t+1} \Vert^2 &- \Vert \tilde \lambda_t \Vert^2 \le 
2 \gamma \big\langle \tilde \lambda_{t+1}, G(x_t) - G(x^*)\big\rangle\\
&- 2 \gamma \big\langle \tilde \lambda_{t+1}, L \tilde \lambda_t \big\rangle
- 2 \gamma \big\langle \tilde \lambda_{t+1}, L \tilde z_t\big\rangle
- \Vert \tilde \lambda_{t+1} - \tilde \lambda_t \Vert^2.
\end{aligned}
\end{equation}
It is clear that
\begin{equation}
\begin{aligned}
\label{z:s:descent}
\Vert \tilde z_{t+1} \Vert^2 - \Vert \tilde z_t\Vert^2 =&~
2\gamma \big\langle \tilde z_t, L \tilde \lambda_t \big\rangle + \gamma^2 \Vert L \tilde \lambda_t \Vert^2.
\end{aligned}
\end{equation}

Combining \eqref{lambda:s:descent} with \eqref{z:s:descent}, we derive
\begin{equation}
\begin{aligned}
\label{property:Vs3}
&V_{s,3}(z_{t+1}, \lambda_{t+1}) - V_{s,3}(z_t, \lambda_t) \\
&\le 2 \gamma \big\langle \tilde \lambda_{t+1}, G(x_t) - G(x^*)\big\rangle
- 2 \gamma \big\langle \tilde \lambda_{t+1}, L \tilde \lambda_t \big\rangle - \Vert \tilde \lambda_{t+1} - \tilde \lambda_t \Vert^2 \\
&~~~ + 3 \gamma^2 \Vert L \tilde \lambda_t \Vert^2
+ 2 \gamma^2 \big\langle \tilde \lambda_{t+1} - \tilde \lambda_t, L^2 \tilde \lambda_t \big\rangle.
\end{aligned}
\end{equation}

Note that $G(\cdot)$ is convex, and $\lambda_t$ is nonnegative for all $t \ge 0$. Plugging \eqref{property:Vs1}, \eqref{property:Vs2} and \eqref{property:Vs3} into \eqref{lyap:V:s}, we obtain
\begin{equation}
\begin{aligned}
\label{property:Vs}
&V_s(x_{t+1}, z_{t+1}, \lambda_{t+1}) - V_s(x_t, z_t, \lambda_t) \\
&\le -2\gamma\mu\|\tilde{x}_t\|^2 -2\gamma \big\langle \tilde{\lambda}_{t+1}, L\tilde{\lambda}_t \big\rangle -\|\tilde{\lambda}_{t+1} - \tilde{\lambda}_t\|^2 - 2\gamma^2\kappa_s\|\nabla G(x_t)^\top\tilde \lambda_t\|^2 \\
&~~~+ 2\gamma(1+\kappa_s) \big\langle \tilde x_t, \nabla G(x_t)^\top (\tilde\lambda_{t+1} - \tilde\lambda_t) \big\rangle +\gamma^2\|\Delta_1(x_t, \lambda_t)\|^2 + 3\gamma^2 \|L\tilde{\lambda}_t\|^2 \\
&~~~  + 2\gamma^2 \big\langle \tilde \lambda_{t+1} - \tilde \lambda_t, L^2\tilde{\lambda}_t \big\rangle +2\gamma\kappa_s \big\langle \tilde{x}_t,
R_t \big\rangle -2\gamma^2\kappa_s \langle \Delta_1(x_t, \lambda^*), \nabla G(x_t)^\top \tilde \lambda_t\rangle \\
&~~~ - 2\gamma^2\kappa_s \big\langle \Delta_1(x_t, \lambda_t), \nabla G(x_t)^\top (\tilde{\lambda}_{t+1} - \tilde{\lambda}_t) + R_t \big\rangle.
\end{aligned}
\end{equation}
In addition, it follows from the Cauchy-Schwarz and Young inequalities that
\begin{align*}
&2\gamma \big\langle \tilde \lambda_{t+1}, L \tilde \lambda_t\big\rangle \ge -\frac 18 \|\tilde \lambda_{t+1}-\tilde \lambda_t\|^2 + \Big(\frac{2\gamma}{k_L} - 8\gamma^2 \Big) \|L \tilde \lambda_t\|^2,\\
&2\gamma(1+\kappa_s)\big\langle \tilde x_t,\nabla G(x_t)^{\top}(\tilde \lambda_{t+1} - \tilde \lambda_t)\rangle
\leq \frac 18 \|\tilde \lambda_{t+1} - \tilde \lambda_t\|^{2} + 8\gamma^{2}(1+\kappa_s)^{2} k_2^{2}\|\tilde x_t\|^{2}, \\
&2\gamma^{2}\langle \tilde \lambda_{t+1} - \tilde \lambda_t, L^{2}\tilde\lambda_t \rangle
\leq \frac 18\|\tilde \lambda_{t+1} - \tilde \lambda_t\|^{2} + 8\gamma^{4} k_L^2 \|L\tilde\lambda_t\|^{2}.
\end{align*}
Since $\kappa_s > 1$, it holds that
\begin{equation}
\begin{aligned}
\label{aux:ieq}
\gamma^{2}\|&\Delta_1(x_t, \lambda_t) \|^{2} - 2\gamma^{2}\kappa_s\|\nabla G(x_t)^\top \tilde \lambda_t\|^{2} - 2\gamma^{2}\kappa_s\langle \Delta_1(x_t, \lambda^*), \nabla G(x_t)^\top \tilde \lambda_t\rangle\\
& =\gamma^{2}\|\Delta_1(x_t, \lambda^*)\|^{2}+2\gamma^{2}(1-\kappa_s)\langle \Delta_1(x_t, \lambda^*), \nabla G(x_t)^\top \tilde \lambda_t \rangle \\
&~~~+ \gamma^{2}(1-2\kappa_s)\|\nabla G(x_t)^\top \tilde \lambda_t\|^{2} \\
&\le -\gamma^{2} \kappa_s\|\nabla G(x_t)^\top \tilde \lambda_t\|^{2} + \gamma^{2} \kappa_s k_5^2 \|\tilde x_t\|^{2}.
\end{aligned}
\end{equation}
Besides, we have
\begin{align*}
&- 2\gamma^{2}\kappa_s \big\langle \Delta_1(x_t, \lambda_t),\nabla G(x_t)^{\top}(\tilde \lambda_{t+1} - \tilde \lambda_t) \big\rangle \\
&~~~~\leq \frac 18 \|\tilde \lambda_{t+1} - \tilde \lambda_t\|^{2} + 16\gamma^{4}\kappa_s^2 k_2^2 k_7^2(\|\tilde x\|^2 + \|\nabla G(x_t)^\top \tilde \lambda_t\|^{2}),\\
&2\gamma\kappa_s \langle \tilde x_t, R_t\rangle
\leq \gamma^{2}\big(2 \kappa_s k_5 k_6 K + \kappa_s^2 k_6^2 K^2\big)\|\tilde x_t\|^{2} + \gamma^{2}\|\nabla G(x_t)^\top \tilde \lambda_t\|^{2}, \\
&- 2\gamma^{2}\kappa_s \big\langle \Delta_1(x_t, \lambda_t), R_t\big\rangle 
\leq 3\gamma^{3} \kappa_s k_6 k_7(k_5 + 1)K (\|\tilde x_t\|^{2} + \|\nabla G(x_t)^\top \tilde \lambda_t\|^{2}).
\end{align*}

Substituting the above estimates into \eqref{property:Vs}, we obtain
\begin{equation}
\begin{aligned}
\label{property:Vs:1}
&V_s(x_{t+1}, z_{t+1}, \lambda_{t+1}) - V_s(x_t, z_t, \lambda_t) \\
&\le - \alpha_1 \|\tilde x_t\|^2  
- \frac 12 \|\tilde \lambda_{t+1}- \tilde \lambda_t\|^2 - \alpha_2 \|\nabla G(x_t)^\top \tilde \lambda_t\|^2 - \alpha_3 \|L \tilde \lambda_t\|^2.
\end{aligned}
\end{equation}
where
\begin{align*}
\alpha_1 =&~ 2\gamma \mu - 8 \gamma^2(1 + \kappa_s)^2 k_2^2 - \gamma^2(2\kappa_s k_5 k_6K + \kappa_s^2 k_6^2 K^2) - \gamma^2 \kappa_s k_5^2 \\
&- 3 \gamma^3 \kappa_s k_6k_7(k_5 + 1)K - 16 \gamma^4 \kappa_s^2 k_2^2 k_7^2, \\
\alpha_2 =&~ \gamma^2(\kappa_s - 1) - 16 \gamma^4 \kappa_s^2 k_2^2 k_7^2 - 3 \gamma^3 \kappa_s k_6k_7(k_5 + 1)K,\\
\alpha_3 =&~ 2 \gamma / k_L - 11\gamma^2 - 8\gamma^4 k_L^2.
\end{align*}
Plugging \eqref{gamma:s:3}-\eqref{gamma:s:5} into \eqref{property:Vs:1}, we derive \eqref{lyap:Vs:onestep}. This completes the proof.
\end{proof}

We are now ready to prove Theorem \ref{thm:semi_global}.

\begin{proof}[Proof of Theorem \ref{thm:semi_global}]
Let $\bar\gamma^s$ be given by \eqref{gamma:s}, and choose any fixed stepsize $\gamma \in (0,\bar\gamma^s)$.
We first show by induction that $(x_t, z_t, \lambda_t) \in \bar{\mathbb B}_{2K_0}(x^*, z^*, \lambda^*)$ for all $t \ge 0$. By the definition of $K_0$, $(x_0,z_0,\lambda_0) \in \bar{\mathbb B}_{K_0}(x^*, z^*, \lambda^*)$.
Suppose that $(x_\tau,z_\tau,\lambda_\tau) \in \bar{\mathbb B}_{2K_0}(x^*, z^*, \lambda^*)$ for $\tau \in \{0, \dots, t\}$. In light of Lemma~\ref{lem:Iter:Bound}, $(x_{\tau + 1},z_{\tau + 1}, \lambda_{\tau + 1}) \in\mathcal K$. It further follows from Lemma \ref{lem:lyap:V:onestep} that 
$V_s(x_{\tau+1},z_{\tau+1}, \lambda_{\tau+1})\le V_s(x_\tau,z_\tau, \lambda_\tau)$ for all $\tau$.
Recalling Lemma \ref{lem:lyap:bound} yields
$$
\frac12\|(x_{t + 1}, z_{t + 1}, \lambda_{t + 1}) - (x^*,z^*, \lambda^*)\|^2 \le V_s(x_{t + 1},z_{t + 1}, \lambda_{t + 1})
\le V_s(x_0, z_0, \lambda_0) \le \frac{3}{2}K_0^2.
$$
Thus, $(x_{t + 1}, z_{t + 1}, \lambda_{t + 1}) \in \bar{\mathbb B}_{2K_0}(x^*,z^*, \lambda^*)$, which completes the induction.

Finally, we show that the trajectory $\{(x_t,z_t,\lambda_t)\}$ converges to the equilibrium set of algorithm \eqref{alg:primal-dual}, which implies the convergence of $\{x_t\}$ to the v-GNE $x^*$.

Define $\mathcal Q_s = \big\{(x, z, \lambda) \in \bar{\mathbb B}_{2K_0}(x^*,z^*, \lambda^*): V_s(x^+, z^+, \lambda^+) =  V_s(x, z, \lambda)\big\}$, where $(x^+,z^+,\lambda^+)$ denotes the one-step update of $(x,z,\lambda)$ under algorithm \eqref{alg:primal-dual}. Let $\mathcal M_s$ be the largest invariant subset of $\bar{\mathcal Q}_s$
By the LaSalle invariance principle \cite{la1976stability}, ${\rm dist}\big((x_t, z_t, \lambda_t), \mathcal M_s\big) \to 0$ as $t \to \infty$.

Let $(\bar x, \bar z, \bar \lambda) \in \mathcal M_s$. It follows from \eqref{lyap:Vs:onestep} that $\bar x = x^*$, $L \bar\lambda = 0$, and $\nabla G(\bar x)^\top (\bar \lambda - \lambda^*) = 0$. Then, $F(\bar x) + \nabla G(\bar x)^\top \bar \lambda = \mathbf{0}$. We further derive $\bar x^+ = \bar x$, and $\bar z^+ = \bar z$. It is also clear that $\bar \lambda^+ = \bar \lambda$ by \eqref{lyap:Vs:onestep}. Therefore, any $(\bar x, \bar z, \bar \lambda) \in \mathcal M_s$ is an equilibrium of \eqref{alg:primal-dual}. Moreover, the trajectory $\{(x_t, z_t, \lambda_t)\}$ approaches $\mathcal M_s$ as $t \to \infty$, and $\{x_t\}$ converges to $x^*$.
This completes the proof.
\end{proof}

%%%%%%
%%%%%%--------------------
\section{Global convergence}
\label{sec:global}

The preceding results indicate that algorithm \eqref{alg:primal-dual} may fail to converge even with an arbitrarily small fixed stepsize. It is therefore natural to consider whether the global convergence can be achieved under a diminishing stepsize sequence. This section is devoted to establishing such a result.

%%%%%%
%%%%%%--------------------
\subsection{Main result}

Consider a variant of algorithm \eqref{alg:primal-dual} with a diminishing stepsize sequence $\{\gamma_t\}$, given by
\begin{equation}
\label{alg:primal-dual:diminish}
\left\{
\begin{aligned}
x_{t+1} &= x_t - \gamma_t F(x_t) - \gamma_t \nabla G(x_t)^\top \lambda_t, \\
z_{t+1} &= z_t + \gamma_t L \lambda_t, \\
\lambda_{t+1} &= P_{\mathbb R^{\tilde m}_+} \big[\lambda_t + \gamma_t G(x_t) - \gamma_t L z_t - \gamma_t L \lambda_t \big].
\end{aligned}
\right.
\end{equation}

We impose the following well-known assumptions.

\begin{assumption}
\label{ass:stepsize}
The stepsize sequence $\{\gamma_t\}$ satisfies
\begin{equation*}
\gamma_t > 0, ~\sum\nolimits_{t} \gamma_t = \infty,~{\rm and}~\sum\nolimits_t \gamma^2_t < \infty.
\end{equation*}
\end{assumption}

\begin{assumption}
\label{ass:g:lips}
For all $i \in \mathcal I$, $g_i$ is $\theta_{g,i}$-Lipschitz continuous.
\end{assumption}

\begin{remark}
We note that Assumption \ref{ass:stepsize} specifies standard stepsize conditions used in stochastic approximation \cite{kushner2003stochastic}. Assumption 5.2 ensures that $\nabla G(x)$ is globally bounded, which plays an important role in the subsequent convergence analysis. Similar assumptions have also been adopted in \cite{nedic2009subgradient, zhu2011distributed, lou2015nash}.
\end{remark}

Even under Assumption \ref{ass:g:lips}, algorithm \eqref{alg:primal-dual} may still fail to converge with a fixed stepsize, as illustrated by the following example.

\begin{example}
\label{Ex:diminishing}
Consider a network with two agents. 
For agent $i \in \{1, 2\}$, the local cost is $J_i(x^i,x^{-i})= (x^i)^2 + \rho x^i x^j$, and the local constraint is
\begin{equation}
g_i(x^i)=
\begin{cases}
\begin{aligned}
&\frac{3}{4}(x^i)^2-\frac{1}{8}(x^i)^4-1, &|x^i|\leq 1, \\
&|x^i|-\frac{11}{8},
& |x^i|\geq 1.
\end{aligned}
\end{cases}
\end{equation}
where $x^i\in\mathbb{R}$, and $0 < \rho < 2$.
The feasible set is $X = \{x \in \mathbb R^2: g_1(x^1) + g_2(x^2) \le 0\}$.
Then the v-GNE is $x^{1,*} = x^{2,*} = 0$.
Take the initial point $x^1_0 = x^2_0 = q_0$, $z_0^1 = z_0^2 = 0$, and $\lambda^1_0 = \lambda^2_0 = \eta_0$. 
By \eqref{alg:primal-dual}, we obtain $x_t^1 = x_t^2 =: q_t$, $\lambda_t^1 = \lambda_t^2 =: \eta_t$, and $z_t^1 = z_t^2 = 0$. Moreover,
\begin{equation}
\begin{cases}
\begin{aligned}
q_{t+1} &= \big(1- (2 + \rho)\gamma\big)q_t- \gamma \nabla g_i(q_t)\eta_t,\\
\eta_{t+1} &= \max\big\{0, \eta_t+\gamma g_i(q_t)\big\}.
\end{aligned}
\end{cases}
\end{equation}
Set $0<\gamma<{2}/{(2 + \rho)}$, and choose $q_0=11/8$ and $\eta_0 = {11(2 - (2 + \rho)\gamma)}/{(8\gamma)}$.
By induction, we obtain  $q_{2t} = 11/8$, $q_{2t + 1} = -11/8$, and $\eta_t = \eta_0$ for all $t$.
Thus, algorithm \eqref{alg:primal-dual} is not convergent.
\end{example}

Example \ref{Ex:diminishing} highlights the limitations of fixed stepsizes. 
We now turn to the global convergence of \eqref{alg:primal-dual:diminish}, which is established in the following theorem.

\begin{theorem}
\label{thm:global:diminish}
Consider algorithm \eqref{alg:primal-dual:diminish}. Let Assumptions \ref{ass:convex}, \ref{ass:monotone}, \ref{ass:graph}, \ref{ass:stepsize} and \ref{ass:g:lips} hold. Then the sequence $\{x_t\}$ converges to the v-GNE of the game \eqref{formulation}.
\end{theorem}

\begin{remark}
Unlike Theorem \ref{thm:semi_global}, Theorem \ref{thm:global:diminish} guarantees the global convergence of \eqref{alg:primal-dual:diminish}, for which the diminishing stepsize sequence is independent of any prescribed compact set of initial conditions. This result is consistent with the asymptotic convergence of its continuous-time counterpart, as the imposed assumptions align with the classical stochastic approximation framework \cite{kushner2003stochastic}.
\end{remark}

\begin{remark}
The convergence of discrete-time primal-dual algorithms for saddle-point problems \cite{nedic2009subgradient}, distributed optimization \cite{zhu2011distributed}, and zero-sum games \cite{lou2015nash} has been investigated. However, these results do not directly establish global convergence of algorithm \eqref{alg:primal-dual:diminish} for distributed GNE seeking in multi-player noncooperative games. Theorem \ref{thm:global:diminish} extends the convergence analysis to this setting.
\end{remark}

%%%%%%
%%%%%%--------------------
\subsection{Convergence analysis}

In this subsection, we prove Theorem \ref{thm:global:diminish}. We begin with the following supporting lemma.

\begin{lemma}
\label{lem:discrete-barbalat}
Let $\{\alpha_t\}_{t=0}^\infty$ and $\{h_t\}_{t=0}^\infty$ be nonnegative sequences satisfying $\sum_{t=0}^{\infty}\alpha_t \\ 
=\infty$, and $\sum_{t=0}^{\infty}\alpha_t h_t<\infty$.
Suppose that there exists a constant $c>0$ such that $|h_{t+1}-h_t| \le c\alpha_t$ for all $t\ge 0$. Then $\lim_{t\to\infty} h_t=0.$
\end{lemma}

\begin{proof}
We first note that $\liminf_{t \to \infty} h_t = 0$, for otherwise $\sum_{t = 0}^\infty \alpha_t h_t = \infty$. It suffices to show $\limsup_{t \to \infty} h_t = 0$.
Suppose, by contradiction, that $\lim\sup_{t \to \infty} h_t > 0$. Then there exists $\varepsilon > 0$ such that $h_t \ge \epsilon$ for infinitely many $t$. Choose $t_1$ such that $h_{t_1} \ge \epsilon$, and recursively define 
$T_k := \min\{t > t_k: h_t < \frac \epsilon 2\}$ and $t_{k+1} = \min\{t > T_k: h_t \ge \epsilon\}$.
Note that $T_k$ must be finite, for otherwise $\sum_{t=0}^{\infty}\alpha_t h_t = \infty$. Then it holds that $h_t \ge \frac \epsilon 2$ for all $t \in [t_k, T_k)$.
Furthermore, we obtain
$$
\frac \epsilon 2 \le h_{t_k} - h_{T_k} \le \sum_{t = t_k}^{T_k - 1} |h_{t+1} - h_t| \le c \sum_{t = t_k}^{T_k - 1} \alpha_t. 
$$

Hence, $\sum_{t = t_k}^{T_k - 1} \alpha_t h_t \ge \frac {\epsilon^2}{4c}$.
Since there are infinitely many disjoint intervals $[t_k, T_k)$, we have $\sum_{k = 1}^\infty \sum_{t = t_k}^{T_k - 1} \alpha_t h_t = \infty$, which contradicts
$\sum_{t=0}^{\infty}\alpha_t h_t < \infty$.
Therefore, $\lim_{t \to \infty} h_t = 0$. The proof is completed.
\end{proof}

We are now ready to prove Theorem \ref{thm:global:diminish}.

\begin{proof}[Proof of Theorem \ref{thm:global:diminish}]
Let $(x^*, z^*, \lambda^*)$ be an equilibrium characterized in Lemma \ref{lem:KKT-EQ}, and define the error variables 
$\tilde x$, $\tilde z$ and $\tilde \lambda$ as in \eqref{error:var}. Construct a Lyapunov function candidate as
\begin{equation}
\label{lyap:V:g}
V_g(x, z, \lambda) = \|\tilde x\|^2 + \|\tilde \lambda\|^2 + \|\tilde z\|^2.
\end{equation} 
Then the proof is divided into the following three steps.

\emph{Step 1:} We investigate the one-step evolution of $V_g$.
Combining \eqref{alg:primal-dual:diminish} with Assumption \ref{ass:monotone} yields
\begin{equation}
\begin{aligned}
\label{property:Vg:1}
\|\tilde x_{t+1}\|^2 &- \|\tilde x_t\|^2
\le -2 \gamma_t \mu \Vert \tilde x_t \Vert^2 
- 2\gamma_t \big\langle \tilde x_t, \nabla G^\top(x_t) \tilde \lambda_t \big\rangle \\
&- 2\gamma_t \big\langle \tilde x_t, \big(\nabla G(x_t) - \nabla G(x^*) \big)^\top \lambda^* \big\rangle + \gamma_t^2 \Vert \Delta_1(x_t, \lambda_t)\Vert^2,
\end{aligned}
\end{equation}
where $\Delta_1(\cdot)$ is defined in \eqref{def:Delta}.

With similar procedures as those of \eqref{lambda:s:descent} and \eqref{z:s:descent}, we obtain
\begin{equation*}
\begin{aligned}
\Vert \tilde \lambda_{t+1} \Vert^2 &- \Vert \tilde \lambda_t \Vert^2 \le 
2 \gamma_t \big\langle \tilde \lambda_{t+1}, G(x_t) - G(x^*)\big\rangle\\
&- 2 \gamma_t \big\langle \tilde \lambda_{t+1}, L \tilde \lambda_t \big\rangle
- 2 \gamma_t \big\langle \tilde \lambda_{t+1}, L \tilde z_t\big\rangle
- \Vert \tilde \lambda_{t+1} - \tilde \lambda_t \Vert^2,\\
\Vert \tilde z_{t+1} \Vert^2 &- \Vert \tilde z_t\Vert^2 =
2\gamma_t \big\langle \tilde z_t, L \tilde \lambda_t \big\rangle + \gamma_t^2 \Vert L \tilde \lambda_t \Vert^2.
\end{aligned}
\end{equation*}
By the convexity of $G$ and the nonnegativity of $\lambda_{t+1}$ and $\lambda^*$, we further derive
\begin{equation*}
\begin{aligned}
\big\langle \tilde \lambda_{t+1}, G(x_t) - G(x^*)\big\rangle \le \big\langle \tilde x_t, \nabla G(x_t)^\top \tilde \lambda_{t+1}\big\rangle + \big\langle \tilde x_t, \big(\nabla G(x_t) - \nabla G(x^*)\big)^\top \lambda^*\big\rangle.
\end{aligned}
\end{equation*}

Combining the above inequalities and recalling \eqref{lyap:V:g}, we obtain
\begin{equation}
\begin{aligned}
\label{property:V:g:1}
&V_g(x_{t+1}, z_{t+1}, \lambda_{t+1}) - V_g(x_t, z_t, \lambda_t) \\
&\le  -2\mu\gamma_t\|\tilde{x}_t\|^2
-\|\tilde \lambda_{t+1} - \tilde \lambda_t\|^2 
+ 2\gamma_t \big\langle\tilde x_t, \nabla G(x_t)^\top(\tilde \lambda_{t+1}-\tilde \lambda_t)\big\rangle \\ 
&~~~ - 2\gamma_t \big\langle \tilde{\lambda}_{t+1}, L\tilde{\lambda}_t \big\rangle
- 2\gamma_t \big\langle \tilde \lambda_{t+1} - \tilde \lambda_t, L\tilde z_t \big\rangle 
+ \gamma_t^2 \Vert \Delta_1(x_t, \lambda_t)\Vert^2
+ \gamma_t^2 \|L\tilde{\lambda}_t\|^2.
\end{aligned}
\end{equation}

Under Assumption \ref{ass:g:lips}, there exists $\theta_g > 0$ such that $\|\nabla G(x)\| \le \theta_g$ for all $x \in \mathbb R^n$. It follows that
$$
\|\Delta_1(x_t, \lambda_t)\|^2 \le 3 \theta^2 \|\tilde x_t\|^2 + 3 \theta_g^2 \|\tilde \lambda_t\|^2
+ 12 \theta_g^2 \|\lambda^*\|^2.
$$
Since $k_L = \|L\|$, it is clear that
\begin{equation*}
\begin{aligned}
2 \gamma_t \big\langle \tilde x_t, \nabla G(x_t)^\top (\tilde \lambda_{t+1} - \tilde \lambda_t)\big\rangle &\le \frac{1}{4}\|\tilde \lambda_{t+1} - \tilde \lambda_t\|^2 + 4\gamma_t^2 \theta_g^2 \|\tilde{x}_t\|^2, \\
- 2\gamma_t \big\langle \tilde\lambda_{t+1}, L \tilde \lambda_t \big\rangle &\le
\frac 14 \|\tilde\lambda_{t+1} - \tilde \lambda_t\|^2 - \Big(\frac{2\gamma_t}{k_L} - 4\gamma_t^2\Big)\|L\tilde \lambda_t\|^2, \\
- 2\gamma_t \big\langle \tilde \lambda_{t+1} - \tilde \lambda_t, L\tilde{z}_t \big\rangle &\le \frac 14 \|\tilde\lambda_{t+1} - \tilde \lambda_t\|^2 + 4\gamma_t^2 k_L^2
\|\tilde{z}_t\|^2.
\end{aligned}
\end{equation*}

Recalling \eqref{property:V:g:1} yields
\begin{equation}
\begin{aligned}
\label{property:V:g:2} 
V_g(x_{t+1}, &z_{t+1}, \lambda_{t+1}) \le (1 + \beta_1 \gamma_t^2) V_g(x_t, z_t, \lambda_t) + \beta_2 \gamma_t^2 \\
&- 2\mu \gamma_t \|\tilde x_t\|^2 - \frac{2 \gamma_t}{k_L} \|L\tilde \lambda_t\|^2 - \frac 1 4 \|\tilde \lambda_{t+1} - \tilde \lambda_t\|^2.
\end{aligned}
\end{equation}
where $\beta_1 = \max\big\{3 \theta^2 + 4 \theta_g^2, 3 \theta_g^2 + 5 k_L^2, 4 k_L^2\big\}$ and
$\beta_2 = 12 \theta_g^2 \|\lambda^*\|^2$.

\emph{Step 2:} We show that $V_g$ is bounded.
By \eqref{property:V:g:2}, we have
$$V_g(x_{t+1}, z_{t+1}, \lambda_{t+1}) \le (1 + \beta_1 \gamma_t^2) V_g(x_t, z_t, \lambda_t) + \beta_2 \gamma_t^2.$$
Resorting to the Gronwall inequality, we further derive
\begin{equation*}
V_g(x_t, z_t, \lambda_t) \le \exp\Bigg(\beta_1 \sum_{k=0}^{t-1}\gamma_k^2\Bigg)\Bigg(V_g(x_0, z_0, \lambda_0) + \beta_2 \sum_{k=0}^{t-1}\gamma_k^2 \Bigg).
\end{equation*}
Under Assumption \ref{ass:stepsize}, $V_g$ is bounded for all $t \ge 0$, and hence, $(x_t, z_t, \lambda_t)$ is  bounded.

\emph{Step 3:} We show that $\{x_t\}$ converges to the v-GNE of \eqref{formulation}. Since $V_g$ is bounded, there exists $\beta_3$ such that $V_g \le \beta_3$ for all $t \ge 0$.
Summing \eqref{property:V:g:2} over $t\ge0$ yields
$$
\sum_{t = 0}^\infty \!\Big(2\mu \gamma_t \|\tilde x_t\|^2 + \frac{2 \gamma_t}{k_L} \|L\tilde \lambda_t\|^2 + \frac 1 4 \|\tilde \lambda_{t+1} - \tilde \lambda_t\|^2\Big) \! \le \! V_g(x_0, z_0, \lambda_0) + (\beta_1 \beta_3 + \beta_2) \sum_{t=0}^\infty \gamma_t^2 \!<\! \infty.
$$
Then $\sum_{t=0}^\infty \gamma_t \|\tilde x_t\|^2 < \infty$. Since the sequence $\{(x_t,z_t,\lambda_t)\}$ is bounded, it follows from \eqref{alg:primal-dual:diminish} that there exists a constant $\beta_4>0$ such that
$|\|\tilde x_{t+1}\|^2 - \|\tilde x_t\|^2| \le \beta_4 \gamma_t$. In light of Lemma \ref{lem:discrete-barbalat}, we 
derive $\lim_{t \to \infty}\|\tilde x_t\| = 0$. Thus, the sequence $\{x_t\}$ converges to $x^*$, and this completes the proof. 
\end{proof}

%%%%%%
%%%%%%--------------------
\section{Convergence under affine constraints}
\label{sec:affine}
The preceding sections have investigated the convergence of \eqref{alg:primal-dual} for games with nonlinear coupling constraints. While a diminishing stepsize sequence ensures global convergence, it may result in a slow convergence rate. This naturally raises the question of whether global convergence can be guaranteed under a fixed stepsize for games with special structures. In this section, we address this question by considering games with affine constraints.

%%%%%%
%%%%%%--------------------
\subsection{Main result}
Consider the game \eqref{formulation} with affine coupling constraints, where each local constraint $g_i$ is of the form
\begin{equation}
\label{affine:constraint}
g_i(x^i) = A_i x^i - b_i, ~\forall i \in \mathcal I.
\end{equation}
where $A_i \in \mathbb R^{m \times n_i}$ and $b_i \in \mathbb R^m$. Accordingly, the feasible set is given by
$X := \{x \in \mathbb R^n \mid \sum_{i \in \mathcal I}(A_i x^i - b_i) \le 0_m\}$.

For convenience, we define $A = {\rm blkdiag}\{A_i\}_{i \in \mathcal I}$ and $b = {\rm col}\{b_i\}_{i \in \mathcal I}$. Under the affine constraints in \eqref{affine:constraint}, algorithm \eqref{alg:primal-dual} is cast into
\begin{equation}
\label{alg:primal-dual:affine}
\left\{
\begin{aligned}
x_{t+1} &= x_t - \gamma F(x_t) - \gamma  A^\top \lambda_t, \\
z_{t+1} &= z_t + \gamma L \lambda_t, \\
\lambda_{t+1} &= P_{\mathbb R^{\tilde m}_+} \big[\lambda_t + \gamma (Ax_t - b) - \gamma L z_t - \gamma L \lambda_t \big],
\end{aligned}
\right.
\end{equation}
where $x(0) \in \mathbb R^n$, $z(0) \in \mathbb R^{\tilde m}$,
and $\lambda(0) \in \mathbb R^{\tilde m}_+$.

\begin{remark}
Affine coupling constraints arise naturally in various applications, including networked Cournot competition \cite{yi2019operator}, network congestion control \cite{Yin2011}, and electric vehicle charging \cite{Paccagnan2019}.
Distributed GNE seeking for \eqref{formulation} with such constraints has been extensively investigated using operator-theoretic approaches, particularly through forward-backward splitting \cite{yi2019operator, pavel2019distributed, pavel2026operator}, forward-backward-forward splitting \cite{Franci2020FBF, belgioioso2022distributed}, and proximal-point methods \cite{Yi2019Proximal}.
In contrast to the aforementioned operator-splitting-based methods, algorithm \eqref{alg:primal-dual} employs explicit primal-dual updates, and does not involve additional forward correction steps or implicit proximal subproblems. Such a simple structure facilitates distributed implementation. However, as far as we know, the convergence of \eqref{alg:primal-dual} has not yet been investigated.
\end{remark}

We now present the main convergence result for \eqref{alg:primal-dual:affine} as follows.

\begin{theorem}
\label{thm:global:affine}
Consider algorithm \eqref{alg:primal-dual:affine} with a fixed stepsize. Let Assumptions \ref{ass:convex}($A_1$) \ref{ass:convex}($A_3$), \ref{ass:monotone} and \ref{ass:graph} hold. Then, there exists $\bar \gamma^a > 0$ such that, for any $\gamma \in (0, \bar \gamma^a)$ and any initial condition $(x_0, z_0, \lambda_0) \in \mathbb R^n \times \mathbb R^{\tilde m} \times \mathbb R_+^{\tilde m}$, the sequence $\{x_t\}$ converges to the v-GNE $x^*$ of game \eqref{formulation} subject to the constraints \eqref{affine:constraint}.
\end{theorem}

\begin{remark}
Theorem \ref{thm:global:affine} indicates that, under a sufficiently small fixed stepsize, the discrete-time primal-dual algorithm is globally convergent to a v-GNE of the game \eqref{formulation} with affine coupling constraints. This strengthens the semiglobal convergence established in Theorem \ref{thm:semi_global}. Although global convergence under fixed stepsizes has been established for several operator splitting approaches \cite{belgioioso2022distributed, pavel2026operator}, these results do not directly apply to algorithm \eqref{alg:primal-dual:affine}. Theorem \ref{thm:global:affine} provides the first global convergence guarantee for this discrete-time primal-dual algorithm, to our knowledge. 
\end{remark}

%%%%%%
%%%%%%--------------------
\subsection{Convergence analysis}
This subsection focuses on the proof of Theorem \ref{thm:global:affine}. Building on the Lyapunov analysis developed in Section \ref{sec:semi}, we exploit the affine structure of the coupling constraints to establish global convergence.

We begin with some preliminary notation. 
Recall that $\theta$ is the Lipschitz constant of $F$ in Assumption \ref{ass:monotone}, and $k_L = \|L\|$.
Let $k_A := \|A\|$, and $k_8 = \max\{\theta, 1\}$. Then
\begin{equation}
\label{ieq:affine}
\forall x \in \mathbb R^n, \forall \lambda \in \mathbb R_+^{\tilde m},~ \|F(x) - F(x^*) + A^\top \tilde\lambda\| \le k_8\big(\|\tilde x\| + \|A^\top \tilde \lambda\|\big).
\end{equation}
We fix a constant $\kappa_a > 1$, and choose the stepsize $\gamma$ such that
\begin{equation}
\label{gamma:a}
\gamma \in (0, \bar \gamma^a), ~\bar \gamma^a:= \min\nolimits_j\{\bar \gamma_j^a\}.
\end{equation}
where
\begin{subequations}
\begin{align}
\bar{\gamma}_1^a &:= \min\Bigg\{\frac{1}{4\kappa_a k_A}, \frac{1}{4 k_L}\Bigg\}, \label{gamma:a:1}\\
\bar{\gamma}_2^a &:= \min\Bigg\{
\frac{\mu}{2\big(8(1+\kappa_a)^2k_A^2+\kappa_a \theta^2\big)},
\sqrt[3]{\frac{\mu}{32 \kappa_a^2 k_A^2k_8^2}}\Bigg\}, \label{gamma:a:2}\\
\bar{\gamma}_3^a &:=\frac{1}{\sqrt{32\kappa_a}\,k_Ak_8}, \label{gamma:a:3}\\
\bar{\gamma}_4^a &:= \min\Bigg\{\frac{1}{22 k_L}, \frac{1}{\sqrt[3]{16}\, k_L}\Bigg\}. \label{gamma:a:4}
\end{align}
\end{subequations}

\begin{remark}
The admissible stepsize $\gamma$ in \eqref{gamma:a} depends on the strong monotonicity and Lipschitz constants of the pseudogradient mapping $F$, the constraint matrix $A$, and the communication topology.
Unlike the bound in \eqref{gamma:s},  the constants involved in \eqref{gamma:a} are independent of the initial conditions. Consequently, a sufficiently small fixed stepsize can be selected to guarantee the convergence of \eqref{alg:primal-dual:affine} from any initial condition, as shown in Theorem \ref{thm:global:affine}. 
\end{remark}

We are now ready to prove Theorem \ref{thm:global:affine}.
\begin{proof}[Proof of Theorem \ref{thm:global:affine}]
Let $(x^*, z^*, \lambda^*)$ be an equilibrium of \eqref{alg:primal-dual:affine}, and define the corresponding error variables $\tilde x$, $\tilde z$ and $\tilde \lambda$ as in \eqref{error:var}.
Following \eqref{lyap:V:s}, we construct a Lyapunov function candidate as
\begin{equation}
\label{lyap:V:a}
V_a(x, z, \lambda) = V_{a,1}(x)  + V_{a,2}(x, \lambda)  + V_{a,3}(z, \lambda),
\end{equation}
where 
\begin{equation*}
\begin{aligned}
V_{a,1}(x) &:= \|\tilde x\|^2, \\
V_{a,2}(x, \lambda) &:= 2\gamma\kappa_a\big\langle \tilde x, A^\top\tilde \lambda\big\rangle,\\
V_{a,3}(z, \lambda) &:= \|\tilde z\|^2 + \|\tilde \lambda\|^2
+ 2\gamma\big\langle \tilde \lambda, L\tilde z\big\rangle.
\end{aligned}
\end{equation*}

By an argument similar to that in Lemma \ref{lem:lyap:bound}, under \eqref{gamma:a:1}, $V_a$ is positive definite and radially unbounded, satisfying
$$
\frac 12 \|(\tilde x, \tilde z, \tilde \lambda)\|^2 \le V_a(x, z, \lambda) 
\le \frac 32  \|(\tilde x, \tilde z, \tilde \lambda)\|^2.
$$

By \eqref{alg:primal-dual:affine}, it is straightforward that
\begin{equation*}
\begin{aligned}
V_{a,1}(x_{t+1}) &- V_{a,1}(x_t) = -2 \gamma \big\langle \tilde x_t, F(x_t) + A^\top \lambda_t \big\rangle + \gamma^2 \Vert F(x_t) + A^\top \lambda_t \Vert^2 \\
&\le -2 \gamma \mu \Vert \tilde x_t \Vert^2 
- 2\gamma \big\langle \tilde x_t, A^\top \tilde \lambda_t \big\rangle 
+ \gamma^2 \Vert F(x_t) + A^\top \lambda_t \Vert^2,
\end{aligned}
\end{equation*}
and moreover,
\begin{equation*}
\begin{aligned}
V_{a,2}(x_{t+1}, \lambda_{t+1}) - V_{a,2}(x_t, \lambda_t)
= 2 \gamma \kappa_a \big\langle \tilde x_t, A^\top(\tilde \lambda_{t+1} 
- \tilde \lambda_t)\big\rangle - 2 \gamma^2 \kappa_a \|A^\top \tilde \lambda_t\|^2 \\
- 2\gamma^2 \kappa_a \big\langle F(x_t) - F(x^*), A^\top \tilde \lambda_t\big\rangle 
- 2 \gamma^2 \kappa_a \big\langle F(x_t) + A^\top \lambda_t, A^\top(\tilde \lambda_{t+1} -\tilde \lambda_t) \big\rangle.
\end{aligned}
\end{equation*}
Following a similar procedure as that of \eqref{property:Vs3}, we obtain
\begin{equation*}
\begin{aligned}
V_{a,3}(z_{t+1}, \lambda_{t+1}) - V_{a,3}(z_t, \lambda_t)
\le 2 \gamma \big\langle \tilde \lambda_{t+1}, A \tilde x_t\big\rangle
- 2 \gamma \big\langle \tilde \lambda_{t+1}, L \tilde \lambda_t \big\rangle& \\
- \Vert \tilde \lambda_{t+1} - \tilde \lambda_t \Vert^2 
+ 3 \gamma^2 \Vert L \tilde \lambda_t \Vert^2
+ 2 \gamma^2 \big\langle \tilde \lambda_{t+1} - \tilde \lambda_t, L^2 \tilde \lambda_t \big\rangle&.
\end{aligned}
\end{equation*}

Substituting the above estimates into \eqref{lyap:V:a}, we derive
\begin{equation}
\begin{aligned}
\label{property:Va:1}
&V_{a}(x_{t+1}, z_{t+1}, \lambda_{t+1}) - V_{a}(x_t, z_t, \lambda_t) \\
&\le -2\gamma\mu\|\tilde x_t\|^2 -2\gamma
\big\langle \tilde \lambda_{t+1}, L\tilde \lambda_t \big\rangle -\|\lambda_{t+1} - \lambda_t\|^2\\
&~~~+ 2\gamma(1+\kappa_a) \big\langle \tilde x, A^\top (\tilde \lambda_{t+1} - \tilde \lambda_t)\big\rangle
+ \gamma^2 \Vert F(x_t) + A^\top \lambda_t \Vert^2 \\
&~~~ -2\gamma^2\kappa_a\|A^\top \tilde \lambda_t\|^2 - 2\gamma^2 \kappa_a \big\langle F(x_t) - F(x^*), A^\top \tilde \lambda_t\big\rangle + 3\gamma^2\|L\tilde \lambda_t \|^2\\
&~~~ - 2\gamma^2\kappa_a \big\langle F(x_t) + A^\top \lambda_t, A^\top (\tilde \lambda_{t+1} - \tilde \lambda_t)\big\rangle + 2\gamma^2 \big\langle \tilde \lambda_{t+1} - \tilde \lambda_t, L^2\tilde \lambda_t \big\rangle.
\end{aligned}
\end{equation}
In addition, we notice that
\begin{equation*}
\begin{aligned}
&2 \gamma \big\langle \tilde \lambda_{t+1}, L \tilde \lambda_t\big\rangle 
\ge - \frac 18 \|\tilde \lambda_{t+1} - \tilde \lambda_t\|^2 + \Big(\frac{2\gamma}{k_L} - 8\gamma^2 \Big) \|L \tilde \lambda_t\|^2, \\
&2\gamma(1+\kappa_a)\big\langle \tilde x_t, A^\top(\tilde \lambda_{t+1} - \tilde \lambda_t)\rangle
\leq \frac 18 \|\tilde \lambda_{t+1} - \tilde \lambda_t\|^{2} + 8\gamma^{2}(1+\kappa_a)^{2} k_A^{2}\|\tilde x_t\|^{2}, \\
&2\gamma^{2} \langle \tilde \lambda_{t+1} - \tilde \lambda_t, L^{2}\tilde\lambda\rangle
\leq \frac 18\|\tilde \lambda_{t+1} - \tilde \lambda_t\|^{2} + 8\gamma^{4}k_L^2 \|L\tilde\lambda_t\|^{2}.
\end{aligned}
\end{equation*}
Similar to \eqref{aux:ieq}, we have
\begin{equation*}
\begin{aligned}
&\gamma^{2}\|F(x_t) + A^\top \lambda_t \|^{2} - 2\gamma^{2}\kappa_a\|A^\top \tilde \lambda_t\|^{2} - 2\gamma^{2}\kappa_a\langle F(x_t) - F(x^*), A^\top \tilde \lambda_t\rangle\\
&\le -\gamma^{2} \kappa_a\|A^\top \tilde \lambda_t\|^{2} + \gamma^{2} \kappa_a \theta^2 \|\tilde x_t\|^{2}.
\end{aligned}
\end{equation*}
Moreover, it is clear that
\begin{equation*}
\begin{aligned}
2\gamma^{2}\kappa_a \big\langle F(x_t) &- F(x^*) + A^\top \tilde \lambda_t, A^\top(\tilde \lambda_{t+1} - \tilde \lambda_t) \big\rangle \\
&\leq \frac 18 \|\tilde \lambda_{t+1} - \tilde \lambda_t\|^{2} + 16\gamma^{4}\kappa_a^2 k_A^2 k_8^2(\|\tilde x\|^2 + \|A^\top \tilde \lambda_t\|^{2}).
\end{aligned}
\end{equation*}

Combining the above estimates with \eqref{property:Va:1}, we obtain
\begin{equation*}
\begin{aligned}
&V_a(x_{t+1}, z_{t+1}, \lambda_{t+1}) - V_a(x_t, z_t, \lambda_t) \\
&\le - \delta_1 \|\tilde x_t\|^2  
- \frac 12 \|\tilde \lambda_{t+1}- \tilde \lambda_t\|^2  - \delta_2 \|A^\top \tilde \lambda_t\|^2 - \delta_3 \|L \tilde \lambda_t\|^2,
\end{aligned}
\end{equation*}
where
\begin{equation*}
\begin{aligned}
\delta_1 &= 2\gamma \mu - 8 \gamma^2(1 + \kappa_a)^2 k_A^2 - \gamma^2\kappa_a \theta^2 - 16 \gamma^4 \kappa_a^2 k_A^2 k_8^2 , \\
\delta_2 &= \gamma^2\kappa_a - 16 \gamma^4 \kappa_a^2 k_A^2 k_8^2, \\
\delta_3 &= 2 \gamma / k_L - 11\gamma^2 - 8\gamma^4 k_L^2.
\end{aligned}
\end{equation*}
Furthermore, plugging \eqref{gamma:a:2}-\eqref{gamma:a:4}, we derive
\begin{equation}
\begin{aligned}
\label{property:Va:3}
&V_a(x_{t+1}, z_{t+1}, \lambda_{t+1}) - V_a(x_t, z_t, \lambda_t) \\
&\le -\gamma \mu\|\tilde x_t\|^2  
- \frac 12 \|\tilde \lambda_{t+1}- \tilde \lambda_t\|^2 
- \frac{\gamma^2\kappa_a}{2} \|\nabla G(x_t)^\top \tilde \lambda_t\|^2 - \frac{\gamma}{k_L}\|L \tilde \lambda_t\|^2.
\end{aligned}
\end{equation}

Since $V_a$ is radially unbounded and nonincreasing along the trajectory of \eqref{alg:primal-dual:affine}, the sequence $\{(x_t,z_t,\lambda_t)\}$ is bounded for all $t \ge 0$.
Let $\mathcal Q_a := \big\{(x, z, \lambda) \in \mathbb R^n \times \mathbb R^{\tilde m} \times \mathbb R^{\tilde m}_+: V_a(x^+, z^+, \lambda^+) = V_a(x, z, \lambda)\big\}$, and $\mathcal M_a$ be the largest invariant subset of $\bar{\mathcal Q}_a$, where $(x^+, z^+, \lambda^+)$ denotes the one-step update of $(x, z, \lambda)$ under \eqref{alg:primal-dual:affine}.
By the LaSalle invariance principle \cite{la1976stability}, $\lim_{t \to \infty}{\rm dist}\big((x_t, z_t, \lambda_t), \mathcal M_a\big) = 0$. Let $(\bar x, \bar z, \bar \lambda) \in \mathcal M_a$. Then 
$\bar x = x^*$ by \eqref{property:Va:3}. Thus, $x_t \to x^*$ as $t \to \infty$, which completes the proof.
\end{proof}

\begin{remark}
For nonlinear coupling constraints considered in Section \ref{sec:semi}, the Lyapunov analysis relies on bounds established over a prescribed compact set, and the admissible stepsize depends on this set. In contrast, under affine constraints, the analysis can be carried out globally using uniform bounds, leading to the global convergence in Theorem \ref{thm:global:affine}.
\end{remark}

%%%%%%
%%%%%%--------------------
\section{Partial-decision information setting}
\label{sec:partial}
We have studied the convergence of \eqref{alg:primal-dual} under a full-decision information setting. 
However, in many practical applications, agents may only communicate with their neighbors and have limited access to the decisions of their opponents \cite{gadjov2018passivity, hu2022distributed}. This motivates the investigation of distributed GNE seeking under a partial-decision information setting, as considered  in \cite{gadjov2018passivity, tatarenko2020geometric, bianchi2021continuous}.

Under this setting, each agent $i$ only receives data from its neighbors via local communication over $\mathcal G_c$, and cannot compute the exact partial-gradient $\nabla_{x^i} J_i(x^i, x^{-i})$ since the complete action profile $x^{-i}$ is not available. To address this limitation, we let each agent maintain local estimates of the other players' actions.
In particular, agent $i$ is endowed with a vector $\mathbf{x}^i = {\rm col}\{\mathbf{x}^{i, j}\}_{j \in \mathcal I} \in \mathbb R^n$, where $\mathbf{x}^{i, i} = x^i$ is agent $i$'s actual strategy, and $\mathbf{x}^{i,j}$ is agent $i$'s estimate of agent $j$'s action.
The estimate components may be different initially, but should reach consensus, i.e., $\lim_{t\to \infty} \Vert \mathbf{x}^i(t) - \mathbf{x}^j(t)\Vert = 0$.
Similar to the full-decision information setting, we introduce a dual variable $\bm \lambda^i \in \mathbb R^m$ and an auxiliary variable $\mathbf{z}^i \in \mathbb R^m$ for agent $i$.
We note that agent $i$ only exchanges $\{\mathbf{x}^i,  \mathbf{z}^i, \bm \lambda^i\}$ with its neighbors over $\mathcal G_c$, enabling a fully distributed implementation.

Since each agent $i$ evaluates its partial gradient $\nabla_{x^i} J_i$ based on its local estimate $\mathbf{x}^i$, we introduce the following extended pseudogradient mapping $\mathbf F$, as a counterpart of $F$ in \eqref{pseudo-gradient}:
\begin{equation}
\label{extend:pseudo-gradient}
\mathbf{F}(\mathbf{x}) := {\rm col}\{\nabla_{x^i} J_i(\mathbf{x}^i)\}_{i \in \mathcal I} \in \mathbb R^n.
\end{equation}
It is clear that $\mathbf{F}(\mathbf{1}_N \otimes x) = F(x)$.
Let $\mathbf{x} := {\rm col}\{\mathbf{x}^i\}_{i \in \mathcal I} \in \mathbb R^{\tilde{n}}$,
$\bm \lambda := {\rm col}\{\bm \lambda^i\}_{i \in \mathcal I}\in \mathbb R^{\tilde{m}}$, $\mathbf{z} := {\rm col}\{\mathbf{z}^i\}_{i \in \mathcal I} \in \mathbb R^{\tilde{m}}$, and $\mathbf{L} = L_G \otimes I_n \in \mathbb R^{\tilde n \times \tilde n}$, where $\tilde n = Nn$. 
To extract agent $i$'s own action $x^i$ from $\mathbf{x}^i$, we define
\begin{equation*}
\begin{aligned}
\mathcal R_i &:= 
\begin{bmatrix}
\mathbf{0}_{n_i \times (n_{< i})}, ~~~~I_{n_i}, ~~~~\mathbf{0}_{n_i \times (n_{> i})}
\end{bmatrix},
\end{aligned}
\end{equation*}
where $n_{< i} = \sum_{j < i, j \in \mathcal I} n_j$, and $n_{> i} = \sum_{j > i, j \in \mathcal I} n_j$.
Furthermore, we let $\mathcal R := {\rm blkdiag}\{\mathcal R_i\}_{i \in \mathcal I}$. Then $\mathcal R \mathbf{x}= x$.

We modify algorithm \eqref{alg:primal-dual} so that each agent $i$ updates its action $x^i$ to minimize its local cost $J_i$, while updating its estimates $\mathbf{x}^{i, -i}$ to achieve consensus with the others.
The resulting fully distributed algorithm is given by
\begin{equation}
\label{alg:primaldual:partial}
\left\{
\begin{aligned}
\mathbf{x}_{t+1} &= \mathbf{x}_t - \gamma_t{\mathcal R}^\top \mathbf{F}(\mathbf{x}_t) - \gamma_t {\mathcal R}^\top \nabla G(\mathcal R \mathbf{x}_t)^\top {\bm \lambda}_t -  \gamma_t c\mathbf{L} \mathbf{x}_t, \\
{\mathbf{z}}_{t+1} &={\mathbf{z}}_t + \gamma_t L {\bm \lambda}_t, \\
{\bm \lambda}_{t+1} &= P_{\mathbb R^{\tilde m}_+} \big[{\bm \lambda}_t + \gamma_t G(\mathcal R \mathbf{x}_t) - \gamma_t L \mathbf{z}_t - \gamma_t L {\bm \lambda}_t \big]
\end{aligned}
\right.
\end{equation}
where $\mathbf{x}(0) \in \mathbb R^{\tilde n}$, $\mathbf{z}(0) \in \mathbb R^{\tilde m}$,
${\bm\lambda}(0) \in \mathbb R^{\tilde m}_+$, $\{\gamma_t\}$ is a positive stepsize sequence and $c > 0$ is the consensus gain.

\begin{remark}
In fact, the continuous-time counterparts of \eqref{alg:primaldual:partial} have been studied in \cite{bianchi2021continuous, li2025passivity}, where asymptotic convergence was established. However, these results do not directly extend to discrete-time settings. To our knowledge, convergence of the discrete-time algorithm \eqref{alg:primaldual:partial} has not yet been investigated in the literature.
\end{remark}

By arguments similar to those in \cite[Lemma 2]{bianchi2021continuous}, the following lemma holds.

\begin{lemma}
\label{KKT:Eq:Partial}
Under Assumptions \ref{ass:convex}, \ref{ass:monotone} and \ref{ass:graph}, $x^*$ is the v-GNE of \eqref{formulation} if and only if there exist $(\mathbf{z}^*, {\bm \lambda}^*)$ such that $(1_N \otimes x^*, \mathbf{z}^*, {\bm \lambda}^*)$ is an equilibrium of \eqref{alg:primaldual:partial}.
\end{lemma}

To facilitate the analysis, we define an augmented pseudogradient mapping as
\begin{equation}
\label{def:F:a}
\mathbf{F}_a(\mathbf{x}) := \mathcal R^\top \mathbf{F}(\mathbf{x}) + c \mathbf{L} \mathbf{x}.
\end{equation}
The next lemma is recalled from \cite[Lemma 3]{bianchi2021continuous}, \cite[Lemma 5]{Meng2026SafeGNE}, and \cite[Lemma 3]{pavel2019distributed}.

\begin{lemma}
\label{lem:partial}
Under Assumptions \ref{ass:convex}, \ref{ass:monotone} and \ref{ass:graph}, the following statements hold.
\begin{enumerate}[label=\roman*)]
\item The mapping $\mathbf{F}$ in \eqref{extend:pseudo-gradient} is $\hat \theta$-Lipschitz continuous for some $\hat \theta \in [\mu, \theta]$.

\item The mapping $\mathbf{F}_a$ in \eqref{def:F:a} is $\bar \theta$-Lipschitz continuous, where $\bar \theta = \hat \theta + c k_L$.

\item Let 
\begin{equation}
\begin{aligned}
\label{def:psi:c}
\Psi = \begin{bmatrix}
\frac{\mu}{N}, & - \frac{\theta + \hat\theta}{2\sqrt{N}} \\
- \frac{\theta + \hat \theta}{2\sqrt{N}}, &c \sigma_2(L_G) - \hat \theta
\end{bmatrix}\!, ~{\rm and}~
c_{\min} = \frac{(\theta + \hat\theta)^2 + 4 \mu \hat\theta}{4 \mu\sigma_2(L_G)},
\end{aligned}
\end{equation}
where $\sigma_2(L_G)$ is the second smallest eigenvalue of the Laplacian matrix $L_G$. Then,
for any $c > c_{\min}$, $\Psi \succ \bm{0}$, and moreover, 
\begin{equation}
\label{restrict:monotone}
\big\langle \mathbf{x} - \mathbf{y}, \mathbf{F}_a(\mathbf{x}) - \mathbf{F}_a(\mathbf{y}) \big\rangle \ge \bar \mu \|\mathbf{x} - \mathbf{y}\|^2,
\end{equation}
where $\bar \mu = \sigma_{\min}(\Psi)$, $\sigma_{\min}$ is the minimal eigenvalue of $\Psi$, $x \in \mathbb R^{\tilde n}$, and $\mathbf{y} = \bm 1_N \otimes y$ for some $y \in \mathbb R^n$
\end{enumerate}
\end{lemma}

Lemma \ref{lem:partial}  establishes the Lipschitz continuity and restricted strong monotonicity of $\mathbf{F}_a$ for $c > c_{\min}$. These properties allow us to extend the convergence analysis in Sections \ref{sec:semi}-\ref{sec:affine} to a partial-decision information setting. 
We summarize the main results in the following theorem.

\begin{theorem}
\label{thm:partial}
Consider algorithm \eqref{alg:primaldual:partial}. Suppose Assumptions \ref{ass:convex}, \ref{ass:monotone} and \ref{ass:graph} hold, and $c > c_{\min}$, where $c_{\min}$ is defined in \eqref{def:psi:c}. Let $x^*$ be the v-GNE of \eqref{formulation} and $\mathbf{x}^* = \bm{1}_N \otimes x^*$. Then the following statements hold.
\begin{enumerate}[label=\roman*)]
\item For every compact set $\bm \Xi_0 \subset \mathbb R^{\tilde n} \times \mathbb R^{\tilde m} \times \mathbb R_+^{\tilde m}$, there exists $\hat \gamma^s > 0$ such that, if $(\mathbf{x}_0, \mathbf{z}_0, \bm\lambda_0) \in \bm \Xi_0$ and $\gamma_t \equiv \gamma \in (0, \hat \gamma^s)$, the sequence $\{\mathbf{x}_t\}$ converges to $\mathbf{x}^*$.

\item If Assumptions \ref{ass:stepsize} and \ref{ass:g:lips} hold, the sequence  $\{\mathbf{x}_t\}$ converges to $\mathbf{x}^*$ for any initial condition $(\mathbf{x}_0, \mathbf{z}_0, \bm\lambda_0) \in \mathbb R^{\tilde n} \times \mathbb R^{\tilde m} \times \mathbb R_+^{\tilde m}$.

\item Suppose that all $g_i$ are affine as in \eqref{affine:constraint}. Then there exists $\hat{\gamma}^a>0$ such that, if $\gamma_t\equiv\gamma\in(0,\hat{\gamma}^a)$,
the sequence $\{\mathbf{x}_t\}$ converges to $\mathbf{x}^*$ for any initial condition
$(\mathbf{x}_0, \mathbf{z}_0, \bm\lambda_0) \in \mathbb R^{\tilde n} \times \mathbb R^{\tilde m} \times \mathbb R_+^{\tilde m}$.
\end{enumerate}
\end{theorem}

\begin{proof}
Define $G_a(\mathbf{x}) := G(\mathcal R \mathbf{x})$. Then, \eqref{alg:primaldual:partial} can be equivalently rewritten as
\begin{equation}
\label{alg:primaldual:partial:rewrite}
\left\{
\begin{aligned}
\mathbf{x}_{t+1} &= \mathbf{x}_t - \gamma_t \mathbf{F}_a(\mathbf{x}_t) - \gamma_t  \nabla G_a(\mathbf{x}_t)^\top {\bm \lambda}_t, \\
{\mathbf{z}}_{t+1} &={\mathbf{z}}_t + \gamma_t L {\bm \lambda}_t, \\
{\bm \lambda}_{t+1} &= P_{\mathbb R^{\tilde m}_+} \big[{\bm \lambda}_t + \gamma_t G_a(\mathbf{x}_t) - \gamma_t L \mathbf{z}_t - \gamma_t L {\bm \lambda}_t \big].
\end{aligned}
\right.
\end{equation}

By Lemma \ref{lem:partial}, $\mathbf{F}_a$ is Lipschitz continuous and restricted strongly monotone. 
In particular, since $\mathbf{x}^* = 1_N\otimes x^*$, it follows from \eqref{restrict:monotone} that
$$
\big\langle \mathbf{x} - \mathbf{x}^*, \mathbf{F}_a(\mathbf{x}) - \mathbf{F}_a(\mathbf{x}^*) \big\rangle \ge \bar \mu \|\mathbf{x} - \mathbf{x}^*\|^2, ~\forall \mathbf{x} \in \mathbb R^{\tilde n}.
$$
This provides the monotonicity estimate required in the Lyapunov analyses of Sections \ref{sec:semi}-\ref{sec:affine}.
Note that $G_a$ inherits the convexity and smoothness properties of $G$.
Consequently, the convergence of \eqref{alg:primaldual:partial} can be analyzed using arguments similar to those in the proofs of Theorems \ref{thm:semi_global}, \ref{thm:global:diminish} and \ref{thm:global:affine}.
The conclusions then follow from Lemma \ref{KKT:Eq:Partial}. This completes the proof.
\end{proof}

\begin{remark}
Theorem \ref{thm:partial} provides convergence guarantees for the discrete-time algorithm \eqref{alg:primaldual:partial} under a partial-decision information setting, complementing the  results in \cite{bianchi2021continuous, li2025passivity}. Unlike the operator splitting approaches in \cite{belgioioso2022distributed, pavel2026operator}, we establish convergence of the explicit primal-dual algorithm without correction steps.
\end{remark}

%%%%%%
%%%%%%--------------------
\section{Simulations}
\label{sec:numerical}

In this section, we carry out numerical simulations to validate the theoretical results.

%%%%%%
%%%%%%--------------------
\subsection{Nonlinear constraints}

Consider a network of $N = 10$ agents. For each agent $i$, the local cost and the constraint are given by $J_i(x^i, x^{-i}) = (x^i)^2 + x^i \sum\nolimits_{j \not= i} x^j$, and
$g_i(x^i) = \sqrt{1 + (x^i)^2} - 2$, where $x^i \in \mathbb R$. The corresponding pseudo-gradient mapping is $F(x) = (I_N + 1_N 1_N^\top) x$, which is $1$-strongly monotone and $(N + 1)$-Lipschitz continuous. Suppose that all agents communicate over a ring graph.

Take the initial point $x^i_0 = q_0$, $\lambda^i_0 = \eta_0$, and $z_0^i = 0$ for all $i$. By \eqref{alg:primal-dual}, we have $x^i_t = q_t$, $\lambda^i_t = \eta_t$, and $z_t^i = 0$ for all $i$. Moreover, it holds that
\begin{equation*}
\left\{
\begin{aligned}
q_{t+1} &= \big(1 - (N + 1) \gamma - \gamma \eta_t/(1 + q_t^2)^{\frac 12}\big) q_t,\\
\eta_{t + 1} &= \max\big\{0, \eta_t + \gamma((1 + q_t^2)^{\frac 12} - 2)\big\}.
\end{aligned}
\right.
\end{equation*}
Let $0 < \gamma < 2/(N + 1)$, and choose $q_0 = \sqrt{3}$ and $\eta_0 = 4/\gamma - 2(N + 1)$. By induction, we obtain
$q_t = (-1)^t \sqrt{3}$ and $\eta_t = \eta_0$ for all $t \ge 0$. Thus, the global convergence of \eqref{alg:primal-dual} cannot be guaranteed, even when the fixed stepsize $\gamma$ is arbitrarily small.

\begin{figure}[htbp]
\centering
\includegraphics[width=0.8\textwidth]{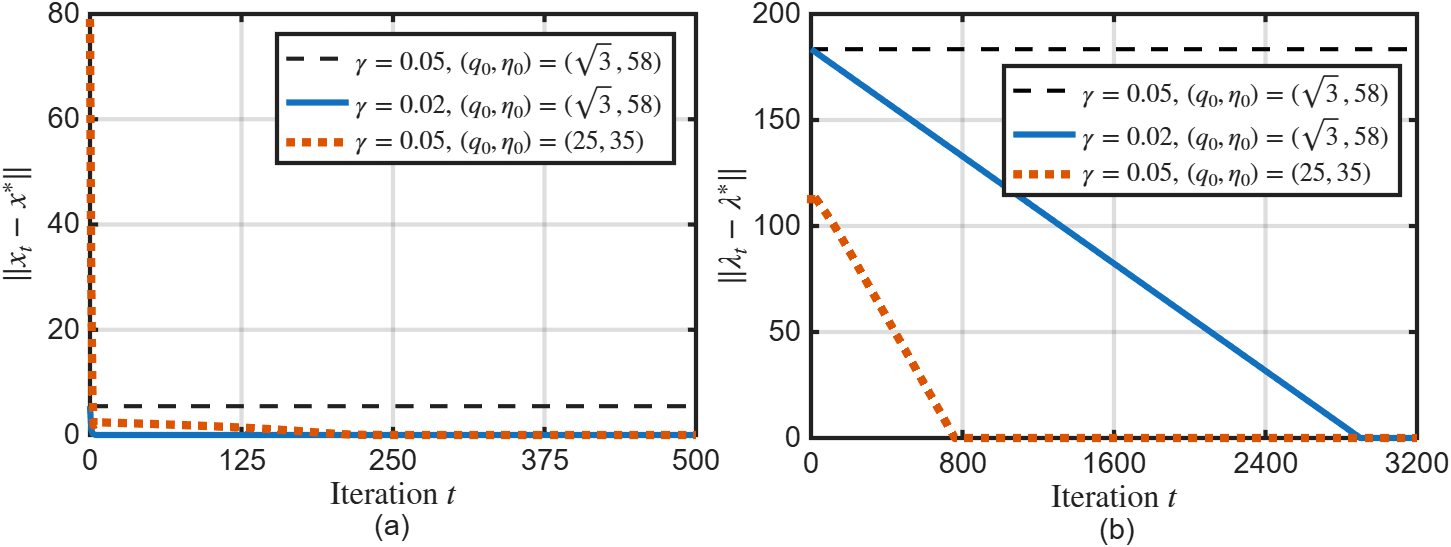}
\caption{Convergence of \eqref{alg:primal-dual} under different stepsizes and initial conditions: (a) $\|x_t - x^*\|$; (b) $\|\lambda_t - \lambda^*\|$.}
\label{fig:1}
\end{figure}

\begin{figure}[htbp]
\centering
\includegraphics[width=0.8\textwidth]{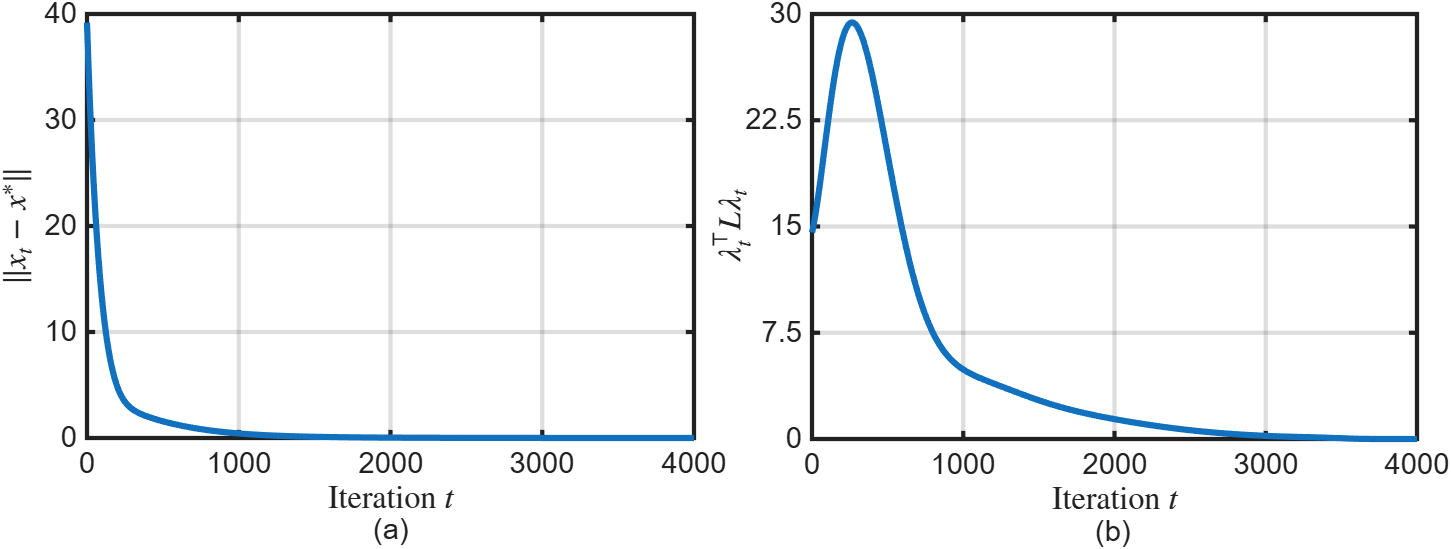}
\caption{Convergence of \eqref{alg:primal-dual}: (a) $\|x_t - x^*\|$; (b) $\lambda_t^\top L \lambda_t$.}
\label{fig:2}
\end{figure}

Fig.~\ref{fig:1} illustrates the convergence of \eqref{alg:primal-dual} under different fixed stepsizes and initial conditions. In particular, for $\gamma=0.05$, the sequence $\{x_t\}$ does not converge when $(q_0, \eta_0) = (\sqrt{3},58)$,  whereas it does when $(q_0,\eta_0)=(25,35)$.
With $(q_0, \eta_0) = (\sqrt{3},58)$, \eqref{alg:primal-dual} is convergent when $\gamma = 0.02$.
These observations highlight the dependence of the admissible stepsize on the initial conditions, consistent with the result in Theorem \ref{thm:semi_global}.
Fig.~\ref{fig:2} shows the convergence of \eqref{alg:primal-dual}, where the agents have different initial values of $x_0^i$, $\lambda_0^i$, and $z_0^i$. It is clear that the sequence $\{x_t\}$ converges to the v-GNE $x^*$, while the local multipliers asymptotically reach consensus.

We now consider a partial-decision information setting, where all agents communicate over a complete graph.
Let $c = 6$, and $0 < \gamma < (31 - \sqrt{895})/33$. Set 
$$\mathbf{x}^{i,i}_0 = \sqrt{3},~\mathbf{x}^{i,j}_0 = - \frac{3\sqrt{3}\gamma}{1 - 3\gamma},~
\bm\lambda_0^i = \frac{4}{\gamma} - \frac{112 - 66 \gamma}{1 - 3\gamma},~\mathbf{z}_0^i = 0,
~\forall i, j \in \mathcal I, i \not= j.$$
By induction, we obtain
$$
\mathbf{x}^{i,i}_t = (-1)^t \sqrt{3},~\mathbf{x}^{i,j}_t= \frac{(-1)^{(t+1)} 3\sqrt{3}\gamma}{1 - 3 \gamma},~
\bm\lambda^i_t = \bm\lambda^i_0,~\mathbf{z}_t^i = 0, ~\forall t \ge 0.
$$
Hence, the sequence $\{x_t\}$ is not convergent in this case. This is consistent with the semi-global convergence result in Theorem \ref{thm:partial}.

\begin{figure}[htbp]
\centering
\includegraphics[width=0.8\textwidth]{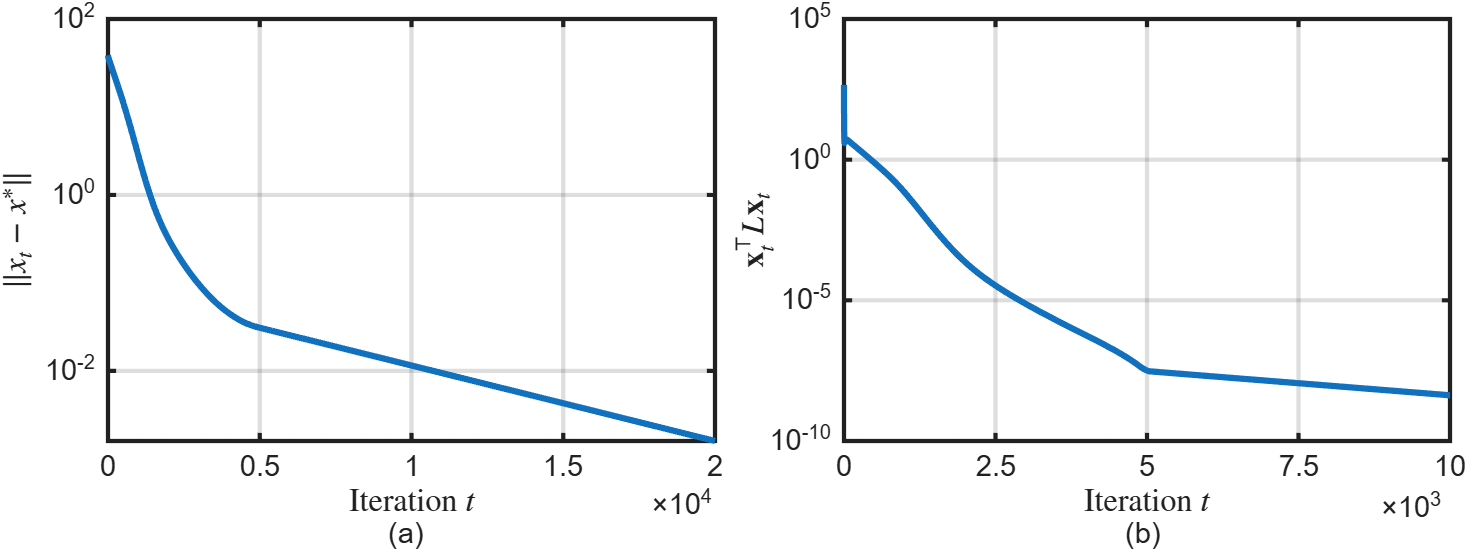}
\caption{Convergence of \eqref{alg:primaldual:partial}: (a) $\|x_t - x^*\|$; (b) $\mathbf{x}_t^\top L \mathbf{x}_t$.}
\label{fig:3}
\end{figure}

Fig.~\ref{fig:3} illustrates the convergence of \eqref{alg:primaldual:partial} over a ring graph for a given initial condition and a fixed stepsize. Obviously, the sequence $\{x_t\}$ converges to $x^*$, and all the local decision estimates reach consensus since $\lim_{t\to \infty}\mathbf{x}^\top_t \bm L \mathbf{x}_t = 0$. 

\begin{figure}[htbp]
\centering
\includegraphics[width=0.8\textwidth]{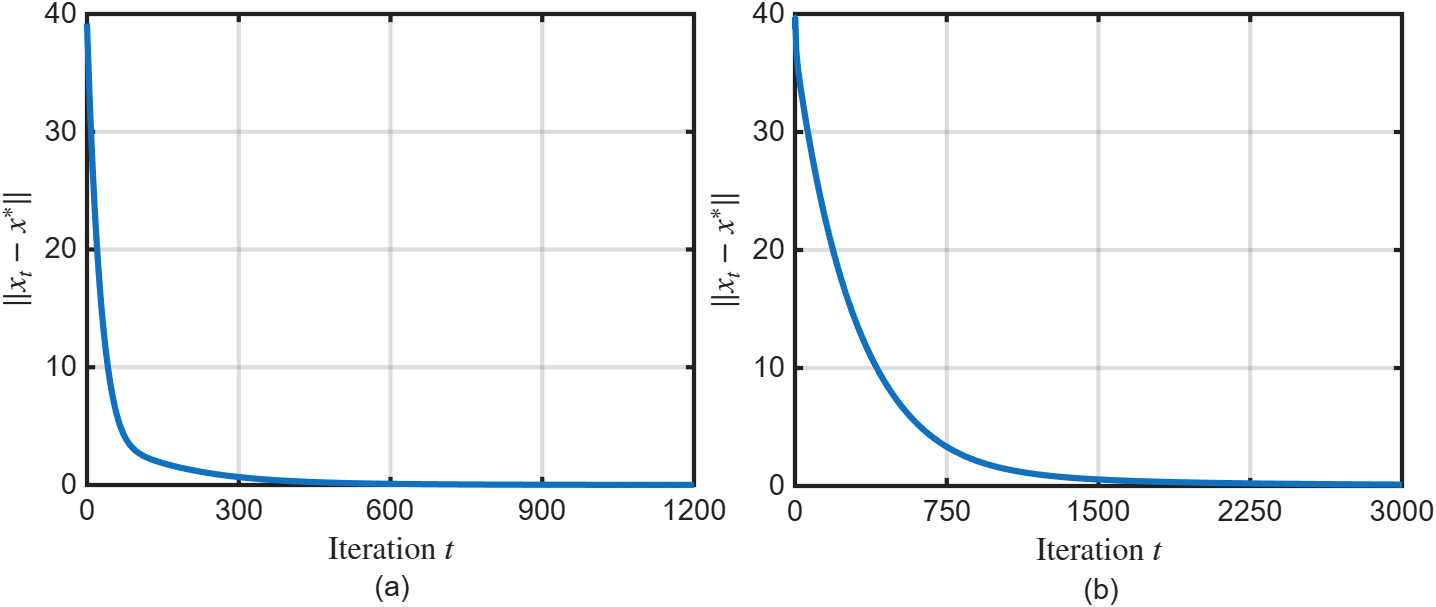}
\caption{Convergence of \eqref{alg:primal-dual:diminish} and \eqref{alg:primaldual:partial} under a diminishing stepsize sequence.}
\label{fig:4}
\end{figure}

Fig.~\ref{fig:4} shows the convergence of \eqref{alg:primal-dual:diminish} and \eqref{alg:primaldual:partial} under a diminishing stepsize sequence, where $\gamma_t = \gamma_0/(1 + \beta t)^\alpha$, $\gamma_0 = 0.003$, $\beta = 0.001$ and $\alpha = 0.6$. Clearly, the sequences $\{x_t\}$ converge to the v-GNE $x^*$ under both full- and partial-decision information settings.

%%%%%%
%%%%%%--------------------
\subsection{Affine constraints}

Consider $N$ firms competing over $m$ markets for a homogeneous commodity.
Firm $i$ produces and delivers $x^i \in \mathbb R^{n_i}$ amount of products to the markets it connect with.
Let $A_i \in \mathbb R^{m \times n_i}$ (with elements $0$ or $1$) be the matrix associated with firm $i$ that specifies the markets it participate in. Define $A := [A_1, \dots, A_N] \in \mathbb R^{m \times n}$ and $x := {\rm col}\{x^i\}$. Then $Ax$ is the total product supply across all markets.
The local cost of firm $i$ is $J_i(x^i, x^{-i}) = c_i(x^i) - P(Ax)^\top A_i x^i$, where $c_i(x^i) = x^{i, \top} Q_i x^i + q_i^\top x^i$ denotes the production cost, $P(Ax) = \bar P - \Xi Ax$ is the price vector,
$\bar P \in \mathbb R^m$ and $\Xi \in \mathbb R^{m \times m}$.
The market capacities impose the affine coupling constraint
$\sum_{i \in \mathcal I} A_i x^i \le \sum_{i \in \mathcal I} r_i$, while the production of firm $i$ is constrained by $0_{n_i} \le x^i \le u_i$, where $r_i \in \mathbb R_+^m$ and $u_i \in \mathbb R_+^{n_i}$.

Let $N = 15$ and $m = 5$. The entries of $\bar P$ are randomly drawn from $(10, 14)$, while 
$\Xi \succ \mathbf{0}$ is chosen with entries in $(1, 2)$. For each $i$, $Q_i \succ \mathbf{0}$ is randomly generated with entries in $(1, 4)$, each component of $q_i$ is randomly from $(0, 2)$, and the entries of $r_i$ and $u_i$ are randomly drawn from $[5, 8]$ and $[6, 14]$.
Each agent is connected to its two nearest neighbors on either side.

\begin{figure}[htbp]
\centering
\includegraphics[width=0.8\textwidth]{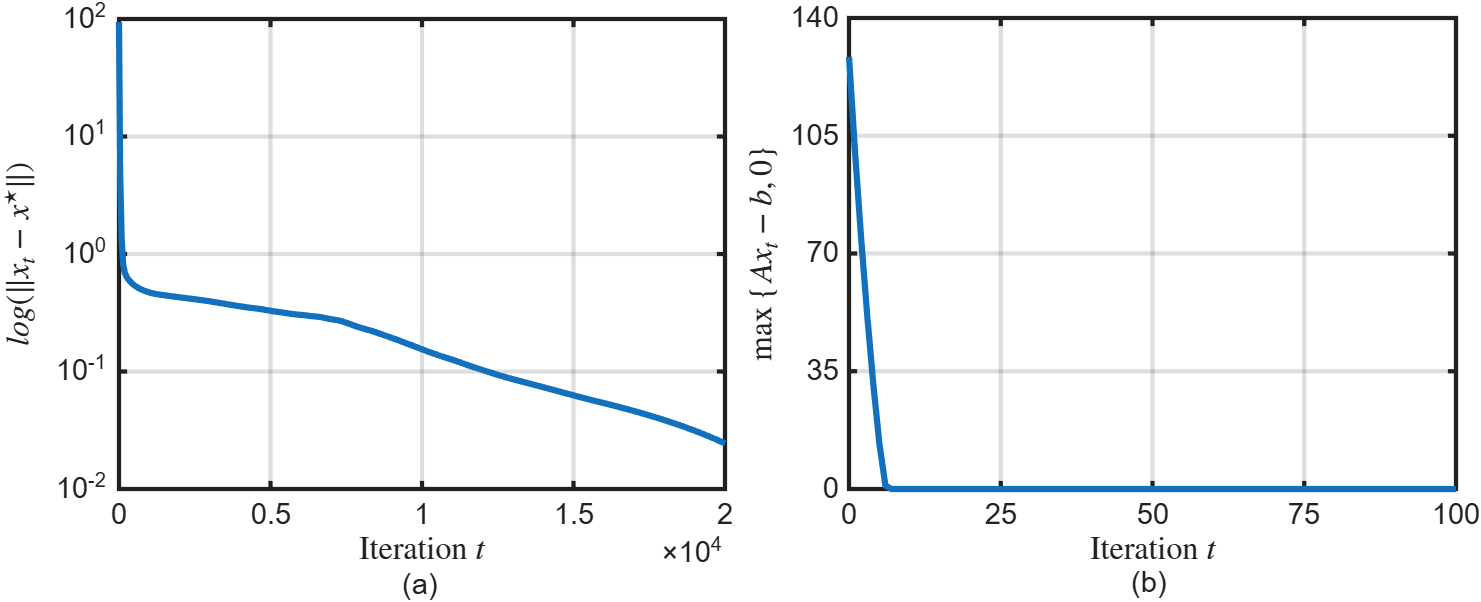}
\caption{Convergence of \eqref{alg:primal-dual:affine} under a full-decision information setting.}
\label{fig:5}
\end{figure}

\begin{figure}[htbp]
\centering
\includegraphics[width=0.8\textwidth]{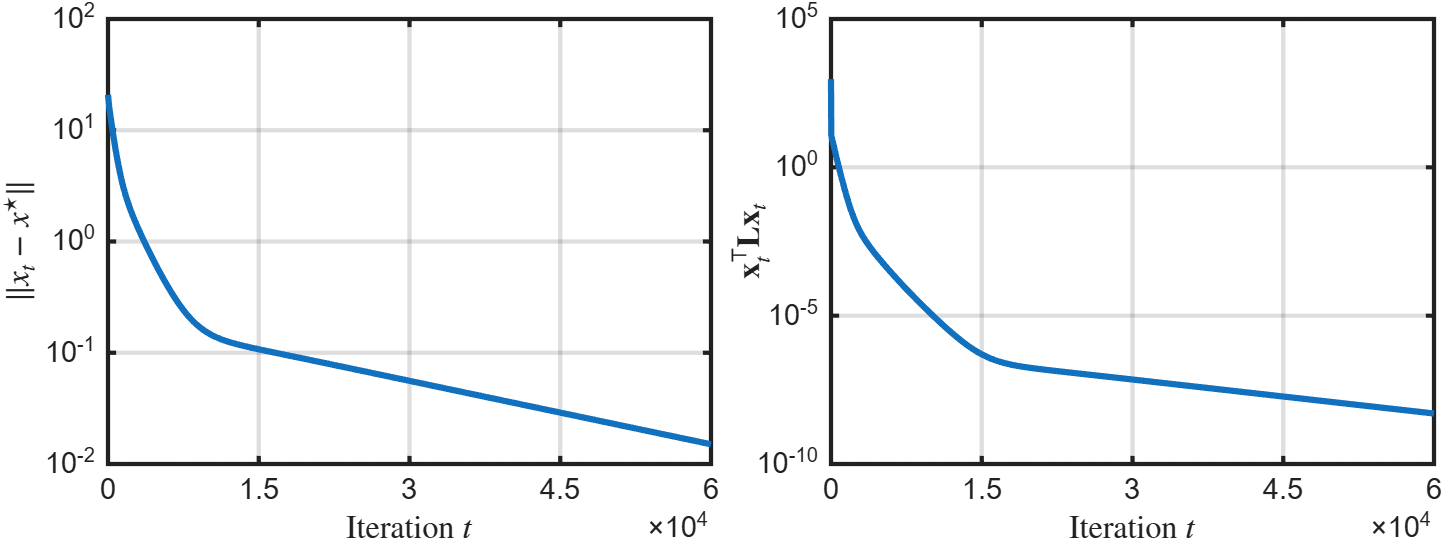}
\caption{Convergence of \eqref{alg:primaldual:partial} under affine constraints.}
\label{fig:6}
\end{figure}

Fig.~\ref{fig:5} shows that, under algorithm \eqref{alg:primal-dual:affine}, $x_t$ converges to the v-GNE $x^*$, while the constraint violation $\max\{Ax_t - b, \mathbf{0}\}$ approaches $\mathbf{0}$.
Fig.~\ref{fig:6} presents the convergence of \eqref{alg:primaldual:partial} with a fixed stepsize under affine constraints. Obviously, $x_t$ converges to $x^*$, while all local estimates $\mathbf{x}^i_t$ reach consensus.
These results validate the theoretical findings in Theorems~\ref{thm:global:affine} and~\ref{thm:partial}.

%%%%%%
%%%%%%--------------------
\section{Conclusion}
\label{sec:conclusion}

This paper established several convergence guarantees of a \allowbreak 
discrete-time distributed primal-dual algorithm for GNE seeking in noncooperative games. For nonlinear coupling constraints, semi-global convergence was shown under a sufficiently small fixed stepsize.
Global convergence was further demonstrated under a stochastic approximation stepsize sequence.
For the special case of affine constraints, global convergence was also guaranteed under a fixed stepsize.
These results were extended to a partial-decision information setting. Finally, numerical simulations were provided for illustration.
Future work will explore globally convergent distributed GNE seeking algorithms with fixed stepsizes under nonlinear coupling constraints.

%\appendix
%\section{An example appendix} 
%\lipsum[71]

% \section*{Acknowledgments}
% We would like to acknowledge the assistance of ....

\bibliographystyle{siamplain}
\bibliography{references}

@article{saad2012game,
	title={Game-theoretic methods for the smart grid: An overview of microgrid systems, demand-side management, and smart grid communications},
	author={Saad, Walid and Han, Zhu and Poor, H Vincent and Basar, Tamer},
	journal={IEEE Signal Processing Magazine},
	volume={29},
	number={5},
	pages={86--105},
	year={2012},
	publisher={IEEE}
}

@article{yin2011nash,
	title={{Nash} equilibrium problems with scaled congestion costs and shared constraints},
	author={Yin, Huibing and Shanbhag, Uday V and Mehta, Prashant G},
	journal={IEEE Transactions on Automatic Control},
	volume={56},
	number={7},
	pages={1702--1708},
	year={2011},
	publisher={IEEE}
}

@article{yan2024incorporation,
	title={Incorporation of likely future actions of agents into pseudo-gradient dynamics of noncooperative games},
	author={Yan, Yuyue and Hayakawa, Tomohisa},
	journal={IEEE Transactions on Automatic Control},
	volume={69},
	number={11},
	pages={7662--7677},
	year={2024},
	publisher={IEEE}
}

@article{johansson2021strategic,
	title={Strategic hub-based platoon coordination under uncertain travel times},
	author={Johansson, Alexander and Nekouei, Ehsan and Johansson, Karl Henrik and M{\aa}rtensson, Jonas},
	journal={IEEE Transactions on Intelligent Transportation Systems},
	volume={23},
	number={7},
	pages={8277--8287},
	year={2021},
	publisher={IEEE}
}

@article{ghaderi2014opinion,
	title={Opinion dynamics in social networks with stubborn agents: equilibrium and convergence rate},
	author={Ghaderi, Javad and Srikant, Rayadurgam},
	journal={Automatica},
	volume={50},
	number={12},
	pages={3209--3215},
	year={2014},
	publisher={Elsevier}
}

@book{facchinei2003finite,
	title={Finite-dimensional variational inequalities and complementarity problems},
	author={Facchinei, Francisco and Pang, Jong-Shi},
	year={2003},
	publisher={Springer}
}

@article{tatarenko2018learning,
	title={Learning generalized {Nash} equilibria in a class of convex games},
	author={Tatarenko, Tatiana and Kamgarpour, Maryam},
	journal={IEEE Transactions on Automatic Control},
	volume={64},
	number={4},
	pages={1426--1439},
	year={2018},
	publisher={IEEE}
}

@article{gao2022continuous,
	title={Continuous-time convergence rates in potential and monotone games},
	author={Gao, Bolin and Pavel, Lacra},
	journal={SIAM Journal on Control and Optimization},
	volume={60},
	number={3},
	pages={1712--1731},
	year={2022},
	publisher={SIAM}
}

@article{arrow1958studies,
	title={Studies in linear and nonlinear programming},
	author={Arrow, Kenneth J and Hurwicz, Leonid and Uzawa, Hirofumi},
	year={1958},
	publisher={Cambridge Univ. Press}
}

@article{polyak1970iterative,
	title={Iterative methods using Lagrange multipliers for solving extremal problems with constraints of the equation type},
	author={Polyak, Boris T},
	journal={USSR Computational Mathematics and Mathematical Physics},
	volume={10},
	number={5},
	pages={42--52},
	year={1970},
	publisher={Elsevier}
}

@article{kallio1999large,
	title={Large-scale convex optimization via saddle point computation},
	author={Kallio, Markku and Rosa, Charles H},
	journal={Operations Research},
	volume={47},
	number={1},
	pages={93--101},
	year={1999},
	publisher={INFORMS}
}

@article{cherukuri2017saddle,
	title={Saddle-point dynamics: conditions for asymptotic stability of saddle points},
	author={Cherukuri, Ashish and Gharesifard, Bahman and Cortes, Jorge},
	journal={SIAM Journal on Control and Optimization},
	volume={55},
	number={1},
	pages={486--511},
	year={2017},
	publisher={SIAM}
}

@article{kia2015distributed,
	title={Distributed convex optimization via continuous-time coordination algorithms with discrete-time communication},
	author={Kia, Solmaz S and Cort{\'e}s, Jorge and Mart{\'\i}nez, Sonia},
	journal={Automatica},
	volume={55},
	pages={254--264},
	year={2015},
	publisher={Elsevier}
}

@article{zeng2018distributed,
	title={Distributed continuous-time algorithms for nonsmooth extended monotropic optimization problems},
	author={Zeng, Xianlin and Yi, Peng and Hong, Yiguang and Xie, Lihua},
	journal={SIAM Journal on Control and Optimization},
	volume={56},
	number={6},
	pages={3973--3993},
	year={2018},
	publisher={SIAM}
}

@article{zeng2019generalized,
	title={Generalized {Nash} equilibrium seeking strategy for distributed nonsmooth multi-cluster game},
	author={Zeng, Xianlin and Chen, Jie and Liang, Shu and Hong, Yiguang},
	journal={Automatica},
	volume={103},
	pages={20--26},
	year={2019},
	publisher={Elsevier}
}

@inproceedings{li2024passivity,
	title={Passivity-based gradient-play dynamics for distributed GNE seeking via parallel feedforward compensation},
	author={Li, Weijian and Pavel, Lacra},
	booktitle={2024 IEEE 63rd Conference on Decision and Control (CDC)},
	pages={2311--2318},
	year={2024},
	organization={IEEE}
}

@article{li2025passivity,
	title={Passivity-based Gradient-Play Dynamics for Distributed Generalized {Nash} Equilibrium Seeking},
	author={Li, Weijian and Pavel, Lacra},
	journal={IEEE Transactions on Automatic Control},
	year={2025},
	publisher={IEEE}
}

@article{pavel2019distributed,
	title={Distributed {GNE} seeking under partial-decision information over networks via a doubly-augmented operator splitting approach},
	author={Pavel, Lacra},
	journal={IEEE Transactions on Automatic Control},
	volume={65},
	number={4},
	pages={1584--1597},
	year={2019},
	publisher={IEEE}
}

@article{facchinei2010generalized,
	title={Generalized {Nash} equilibrium problems},
	author={Facchinei, Francisco and Kanzow, Christian},
	journal={Annals of Operations Research},
	volume={175},
	number={1},
	pages={177--211},
	year={2010},
	publisher={Springer}
}

@article{kulkarni2012variational,
	title={On the variational equilibrium as a refinement of the generalized {Nash} equilibrium},
	author={Kulkarni, Ankur A and Shanbhag, Uday V},
	journal={Automatica},
	volume={48},
	number={1},
	pages={45--55},
	year={2012},
	publisher={Elsevier}
}

@article{hu2022distributed,
	title={Distributed {Nash} equilibrium seeking: Continuous-time control-theoretic approaches},
	author={Hu, Guoqiang and Pang, Yipeng and Sun, Chao and Hong, Yiguang},
	journal={IEEE Control Systems Magazine},
	volume={42},
	number={4},
	pages={68--86},
	year={2022},
	publisher={IEEE}
}

@article{yi2019operator,
	title={An operator splitting approach for distributed generalized {Nash} equilibria computation},
	author={Yi, Peng and Pavel, Lacra},
	journal={Automatica},
	volume={102},
	pages={111--121},
	year={2019},
	publisher={Elsevier}
}

@article{bianchi2021continuous,
	title={Continuous-time fully distributed generalized {Nash} equilibrium seeking for multi-integrator agents},
	author={Bianchi, Mattia and Grammatico, Sergio},
	journal={Automatica},
	volume={129},
	pages={109660},
	year={2021},
	publisher={Elsevier}
}

@article{gadjov2018passivity,
	title={A passivity-based approach to {Nash} equilibrium seeking over networks},
	author={Gadjov, Dian and Pavel, Lacra},
	journal={IEEE Transactions on Automatic Control},
	volume={64},
	number={3},
	pages={1077--1092},
	year={2018},
	publisher={IEEE}
}

@article{tatarenko2020geometric,
	title={Geometric convergence of gradient play algorithms for distributed {Nash} equilibrium seeking},
	author={Tatarenko, Tatiana and Shi, Wei and Nedi{\'c}, Angelia},
	journal={IEEE Transactions on Automatic Control},
	volume={66},
	number={11},
	pages={5342--5353},
	year={2020},
	publisher={IEEE}
}

@book{la1976stability,
	title={The stability of dynamical systems},
	author={La Salle, Joseph P},
	year={1976},
	publisher={SIAM}
}

@article{nedic2009subgradient,
	title={Subgradient methods for saddle-point problems},
	author={Nedi{\'c}, Angelia and Ozdaglar, Asuman},
	journal={Journal of optimization theory and applications},
	volume={142},
	number={1},
	pages={205--228},
	year={2009},
	publisher={Springer}
}

@inproceedings{wang2010control,
	title={Control approach to distributed optimization},
	author={Wang, Jing and Elia, Nicola},
	booktitle={2010 48th Annual Allerton Conference on Communication, Control, and Computing (Allerton)},
	pages={557--561},
	year={2010},
	organization={IEEE}
}

@article{lei2016primal,
	title={Primal--dual algorithm for distributed constrained optimization},
	author={Lei, Jinlong and Chen, Han-Fu and Fang, Hai-Tao},
	journal={Systems \& Control Letters},
	volume={96},
	pages={110--117},
	year={2016},
	publisher={Elsevier}
}

@article{gharesifard2013distributed,
	title={Distributed convergence to {Nash} equilibria in two-network zero-sum games},
	author={Gharesifard, Bahman and Cort{\'e}s, Jorge},
	journal={Automatica},
	volume={49},
	number={6},
	pages={1683--1692},
	year={2013},
	publisher={Elsevier}
}

@article{lou2015nash,
	title={Nash equilibrium computation in subnetwork zero-sum games with switching communications},
	author={Lou, Youcheng and Hong, Yiguang and Xie, Lihua and Shi, Guodong and Johansson, Karl Henrik},
	journal={IEEE Transactions on Automatic Control},
	volume={61},
	number={10},
	pages={2920--2935},
	year={2015},
	publisher={IEEE}
}

@article{zhu2011distributed,
	title={On distributed convex optimization under inequality and equality constraints},
	author={Zhu, Minghui and Martinez, Sonia},
	journal={IEEE Transactions on Automatic Control},
	volume={57},
	number={1},
	pages={151--164},
	year={2011},
	publisher={IEEE}
}

@book{kushner2003stochastic,
	title={Stochastic approximation and recursive algorithms and applications},
	author={Kushner, Harold J and Yin, G George},
	year={2003},
	publisher={Springer}
}

@article{pavel2026operator,
	title={On operator theory and applications in game theory},
	author={Pavel, Lacra},
	journal={Annual Review of Control, Robotics, and Autonomous Systems},
	volume={9},
	year={2026},
	publisher={Annual Reviews}
}

@article{belgioioso2022distributed,
	title={Distributed generalized {Nash} equilibrium seeking: An operator-theoretic perspective},
	author={Belgioioso, Giuseppe and Yi, Peng and Grammatico, Sergio and Pavel, Lacra},
	journal={IEEE Control Systems Magazine},
	volume={42},
	number={4},
	pages={87--102},
	year={2022},
	publisher={IEEE}
}

@article{Meng2026SafeGNE,
	author={Yihan Meng and Weijian Li and Lacra Pavel},
	title={Safe Distributed Generalized {Nash} Equilibrium Seeking via Control Barrier Functions},
	journal={arXiv preprint arXiv:2608.24077},
	year={2026}
}

@article{bin2024semiglobal,
	title={Semiglobal exponential stability of the discrete-time {Arrow-Hurwicz-Uzawa} primal-dual algorithm for constrained optimization},
	author={Bin, Michelangelo and Notarnicola, Ivano and Parisini, Thomas},
	journal={Mathematical Programming},
	volume={208},
	number={1},
	pages={629--660},
	year={2024},
	publisher={Springer}
}

@article{ren2026exponentially,
	title={An exponentially stable discrete-time primal-dual algorithm for distributed constrained optimization},
	author={Ren, Xiaoxing and Bin, Michelangelo and Notarnicola, Ivano and Parisini, Thomas},
	journal={IEEE Transactions on Automatic Control},
	year={2026},
	publisher={IEEE}
}

@book{grne2013nonlinear,
	title={Nonlinear model predictive control: theory and algorithms},
	author={Grne, Lars and Pannek, Jrgen},
	year={2013},
	publisher={Springer Publishing Company, Incorporated}
}

@article{Paccagnan2019,
	author  = {Dario Paccagnan and Basilio Gentile and
	Francesca Parise and Maryam Kamgarpour and
	John Lygeros},
	title   = {{Nash} and Wardrop Equilibria in Aggregative Games
	With Coupling Constraints},
	journal = {IEEE Transactions on Automatic Control},
	volume  = {64},
	number  = {4},
	pages   = {1373--1388},
	year    = {2019}
}

@article{Yin2011,
	author  = {Huibing Yin and Uday V. Shanbhag and Prashant G. Mehta},
	title   = {{Nash} Equilibrium Problems With Scaled Congestion Costs and Shared Constraints},
	journal = {IEEE Transactions on Automatic Control},
	volume  = {56},
	number  = {7},
	pages   = {1702--1708},
	year    = {2011}
}

@article{Yi2019Proximal,
	author  = {Peng Yi and Lacra Pavel},
	title   = {Distributed Generalized {Nash} Equilibria Computation of Monotone Games via Double-Layer Preconditioned Proximal-Point Algorithms},
	journal = {IEEE Transactions on Control of Network Systems},
	volume  = {6},
	number  = {1},
	pages   = {299--311},
	year    = {2019}
}

@inproceedings{Franci2020FBF,
	author    = {Barbara Franci and Mathias Staudigl and Sergio Grammatico},
	title     = {Distributed Forward-Backward (Half) Forward Algorithms for Generalized {Nash} Equilibrium Seeking},
	booktitle = {2020 European Control Conference (ECC)},
	pages     = {1274--1279},
	year      = {2020},
	publisher = {IEEE},
}
\end{document}